\documentclass{amsart}
\usepackage{amsaddr}
\usepackage[T1]{fontenc}
\usepackage{amsmath,amssymb,amscd,amsthm}
\usepackage[latin1]{inputenc}
\usepackage{latexsym}
\usepackage[dvips]{graphicx}
\usepackage{epsfig}
\usepackage{graphics}
\usepackage{colortbl}
\usepackage{color}
\usepackage{latexsym}
\usepackage{ae}
\usepackage[cm]{aeguill}

\begin{document}
\setcounter{secnumdepth}{2}
\setcounter{tocdepth}{2}
\newtheorem{theorem}{Theorem}[section]
\newtheorem{definition}[theorem]{Definition}
\newtheorem{lemma}[theorem]{Lemma}
\newtheorem{proposition}[theorem]{Proposition}
\newtheorem{question}[theorem]{Question}
\newtheorem{conjecture}[theorem]{Conjecture}
\newtheorem{corollary}[theorem]{Corollary}
\newtheorem{gstatement}{Question}
\newcommand{\mf}{\mathfrak}
\newcommand{\mb}{\mathbb}
\newcommand{\ol}{\overline}
\newcommand{\la}{\langle}
\newcommand{\ra}{\rangle}
\newcommand{\cotanh}{\mathrm{cotanh}}

\newtheorem{Alphatheorem}{Theorem}
\renewcommand{\theAlphatheorem}{\Alph{Alphatheorem}}

\newcommand{\EM}{\ensuremath}
\newcommand{\norm}[1]{\EM{\left\| #1 \right\|}}

\newcommand{\pa}{\parallel}

\newcommand{\I}{\textsc{I}}
\newcommand{\II}{\textsc{I\hspace{-0.05 cm}I}}
\newcommand{\III}{\textsc{I\hspace{-0.05 cm}I\hspace{-0.05 cm}I}}

\title{Fuchsian polyhedra in Lorentzian space-forms}

\author{ February 2009}

\author{Fran\c{c}ois Fillastre }
\address{Universit\'e de Cergy-Pontoise \\ UMR CNRS 8088, d\'epartement de math\'ematiques, \\
F-95000 Cergy-Pontoise \\
FRANCE \\
{\tt francois.fillastre@u-cergy.fr} \\
{\tt http://www.u-cergy.fr/ffillast/}
}

\maketitle

\begin{abstract}

Let $S$ be   a compact surface of genus $>1$, and $g$ be a  metric on  $S$ of constant curvature $K\in\{-1,0,1\}$ with   conical singularities of negative singular curvature. When $K=1$ we add the condition that the lengths of the contractible geodesics are $>2\pi$. 
 We prove that there exists a convex polyhedral surface $P$ in the Lorentzian space-form of curvature $K$  and a group $G$ of isometries of this space such that the induced metric on the quotient $P/G$ is isometric to $(S,g)$. Moreover, the pair $(P,G)$ is unique (up to global isometries) among a particular class of convex polyhedra, namely Fuchsian polyhedra. 
 
This extends  theorems of  A.D. Alexandrov and Rivin--Hodgson \cite{alex42,RivHod} concerning the sphere to the higher genus cases, and it is also the polyhedral version of a theorem of Labourie--Schlenker \cite{SchLab}.
\end{abstract}



\section{Definitions and statements}

In Subsection  \ref{subsecDefs} we recall standard definitions and results. The main result and original definitions of the paper are stated 
in Subsection \ref{subsec:FuchsPol}.

\subsection{Metrics with conical singularities and convex polyhedra.}\label{subsecDefs}

Let $M_K^-$ be the Lorentzian space-form of dimension $3$ with constant
curvature $K$, $K\in\{-1,0,1\}$:  $M_0^-$ is the Minkowski space, $M_1^-$ is the de Sitter space and $M_{-1}^-$ is the anti-de Sitter space. 
A \emph{convex polyhedron} is an intersection of half-spaces of $M_K^-$. The number of half-spaces may be infinite, but the intersection is asked to be  locally finite:  each face must be a polygon with a finite number of vertices, and  the number of edges at each vertex must be finite. A \emph{polyhedron} is a connected union of convex polyhedra. 
 A \emph{polyhedral surface} is the boundary of a polyhedron and a \emph{convex polyhedral surface} is the boundary of a convex polyhedron.
A \emph{convex (polyhedral) cone}
in $M_K^-$ is a  convex polyhedral surface with only one vertex. In the de Sitter case, we will call \emph{polyhedral surface of hyperbolic type} a polyhedral surface dual of a hyperbolic polyhedral surface (the definition of duality is recalled in Section \ref{section backgrounds}).
The sum of the angles between the edges on the faces of a space-like convex cone in $M_K^-$ (of hyperbolic type for $M_1^-$) is strictly
greater than $2\pi$.

A  \emph{metric of curvature $K$ with conical singularities of negative singular curvature} on a compact surface $S$ is a
(Riemannian) metric of constant curvature $K$ on $S$ minus $n$
points $(x_1,\ldots,x_n)$ such that the neighbourhood of each $x_i$
is isometric to the induced metric on the neighbourhood of the vertex of a convex
cone in $M_K^-$. The $x_i$ are called the \emph{singular points}. By
definition the set of singular points is discrete, hence finite since
the surface is compact. The angle $\alpha_i$ around a singular point $x_i$ is the \emph{cone-angle} at this point and the value $(2\pi-\alpha_i)$ is the \emph{singular curvature} at $x_i$.

Let $P$ be a convex polyhedral surface in $M_1^-$ of hyperbolic type homeomorphic to the sphere (note that the de Sitter space is not contractible). The
induced metric on $P$ is isometric to a   spherical  metric with conical singularities of negative singular curvature on
the sphere. Moreover the lengths of the closed geodesics for this metric are $>2\pi$, see  \cite{RivHod}. The  following theorem says that all the metrics of this kind can be obtained by such polyhedral surfaces: 
\begin{theorem}[{Rivin--Hodgson, \cite{RivHod}}]\label{thmrivhod}
Each  spherical metric on the
sphere with conical singularities of negative singular curvature such that the lengths of its closed geodesics are $>2\pi$   can be isometrically embedded in the de Sitter space as  a unique (up to global isometries) convex  polyhedral surface of hyperbolic type  homeomorphic to the sphere. 
\end{theorem}

This result was extended to the cases where $P$ is not of hyperbolic type in \cite{schhilb,schlorentz} (actually the results contained in these references cover larger classes of metrics on the sphere).
It is an extension to negative singular curvature of a famous theorem of A.D. Alexandrov. We denote by $M_K^+$ the Riemannian space-form of dimension $3$ with constant curvature $K$. A conical singularity with \emph{positive singular curvature} is a point which has a neighbourhood  isometric to the induced metric on the neighbourhood of the vertex of a convex
cone in $M_K^+$. The sum of the angles between the edges on the faces of a convex cone in $M_K^+$ is strictly
between $0$ and $2\pi$.
\begin{theorem}[{A.D. Alexandrov, \cite{alex42,Aleks}}]\label{thmalex}
Each metric of curvature $K$ on the sphere  with conical singularities of positive singular curvature can be isometrically embedded in $M_K^+$ as (the boundary of) a unique (up to global isometries) convex compact polyhedron.
\end{theorem}

By the Gauss--Bonnet formula, we know that there doesn't exist other constant curvature metrics with conical singularities of constant sign singular curvature on the sphere than the ones described in Theorem  \ref{thmrivhod} and Theorem \ref{thmalex}. In the present paper we extend these results to the cases of surfaces of higher genus (actually $>1$).

\subsection{Fuchsian polyhedra.}\label{subsec:FuchsPol}

An \emph{invariant polyhedral surface} is a pair $(P,F)$, where $P$ is a polyhedral surface in $M_K^-$ and 
$F$ is a discrete group of isometries of $M_K^-$ such that $F(P)=P$ and $F$ acts freely on $P$. The group $F$ is called 
the \emph{acting group}.
If there exists an invariant polyhedral surface $(P,F)$ in $M_K^-$ such that the induced metric on $P/F$ is isometric to a metric $h$ on a surface $S$,  we say that $P$ \emph{realises} the metric $h$ (the singular points of $h$ will correspond to the vertices of
$P$, and $F$ will be isomorphic to the fundamental group of $S$).  For example, if $S$ is the sphere the acting group is the trivial one.

We denote by $\mathrm{Isom}^+_+(M_K^-)$ the group of orientation preserving and time orientation preserving 
isometries of $M_K^-$.
We call a \emph{Fuchsian group of $M_K^-$}  a  discrete subgroup of $\mathrm{Isom}^+_+(M_K^-)$ fixing a point $c_K$ and acting freely cocompactly on the time-like units vectors at $c_K$.
Such a group also leaves invariant and acts freely cocompactly on all the surfaces in the future-cone of $c_K$ which are at constant distance from $c_K$. These surfaces have the properties to be space-like, strictly convex, umbilical and complete. 
The induced metric on exactly one of these surfaces, denoted  by  $O_K$, is hyperbolic.
A \emph{Fuchsian polyhedron} of $M_K^-$ is an invariant  polyhedral surface $(P,F)$ where $P$ is a polyhedral surface contained in the future-cone of $c_K$ and $F$ is a Fuchsian group of $M_K^-$. In this paper we assume moreover that a Fuchsian polyhedron is space-like. Example of Fuchsian polyhedra are given in the proof of Lemma~\ref{lem:exFuchs}.

\begin{Alphatheorem}\label{realisation lorentz}
Let $S$ be a compact surface of genus $>1$.
\begin{itemize}
\item[1)] A  spherical metric with conical singularities of negative singular curvature on $S$ such that the lengths  of its 
closed contractible geodesics are  $>2\pi$ is realised by a unique (up to global isometries)  convex Fuchsian
polyhedron in the de Sitter space.

\item[2)] A  flat metric with conical singularities of negative singular curvature on $S$ is realised by a unique  (up to global isometries)  convex Fuchsian
polyhedron in the  Minkowski space.

\item[3)] A  hyperbolic  metric with conical singularities of negative singular curvature on $S$ is realised by a unique  (up to global isometries)  convex Fuchsian
polyhedron in the  anti-de Sitter space.
\end{itemize}
\end{Alphatheorem}

In the statements above, uniqueness must be understood as the uniqueness among the class of convex polyhedra described in the statements. Otherwise the statements are false. It is clear if the polyhedra are not required to be convex, and it should be 
easy to construct other examples of invariant (convex) polyhedra realising the described metrics on $S$.
The necessity of the condition on the lengths of the geodesics in the spherical case will be explained in Section \ref{section backgrounds}.

The part 1) of Theorem~\ref{realisation lorentz} was already done in \cite[thm 4.22]{Schpoly}. The other parts 
are new from what I know. The part 2) of Theorem~\ref{realisation lorentz} answers a question asked in \cite{Schpolygone}.  
The analog of Theorem~\ref{realisation lorentz} for smooth metrics was proved in \cite{SchLab}.

\subsection{Outline of the proof --- organisation of the paper.}

Actually the general outline of the proof is very classical,  starting from Alexandrov's work, and very close to the one used in \cite{artrealisationhyp}. 
Roughly speaking, the idea is to endow with   suitable topology both the space of convex Fuchsian polyhedra of $M_K^-$ and the space of corresponding  metrics, and to show that the map from one to the other given by the induced metric is a homeomorphism.

The difficult steps are (always) to show  local injectivity and  properness of the maps ``induced
metric''. The local injectivity is equivalent to  statements about infinitesimal rigidity of convex
Fuchsian polyhedra.  The Minkowski case is already known and the others cases will follow by the use of  the so-called infinitesimal Pogorelov maps.

In the remainder of this section we present some consequences and possible extensions of Theorem~\ref{realisation lorentz} --- note that the theorems proved in this paper are labeled with letters instead of numbers.

In the following section we present ``projective models'' of the space-forms and describe the shape of convex Fuchsian polyhedra. We also recall some facts about Teichm\"uller space which will be used in the sequel.
Section~\ref{rigidite lorentz} is devoted to the proof of the infinitesimal rigidity of convex Fuchsian polyhedra among convex Fuchsian polyhedra.
Subsection~\ref{Set of Fuchsian polyhedra}  studies the topologies of the spaces of polyhedra, Subsection~\ref{Set of metrics}  studies the ones of the spaces of metrics,  and Subsection~\ref{final steps} gives a sketch of the proof of Theorem~\ref{realisation lorentz}.
Finally, we prove in Section~\ref{properness} the properness of the maps ``induced metric'', that was the last thing to check according to the sketch of the preceding section.

\subsection{Fuchsian polyhedra in hyperbolic space --- Towards a general result.}

  In the hyperbolic space a Fuchsian polyhedron is a polyhedral surface invariant under the action of a Fuchsian group of hyperbolic space $\mathbb{H}^3$. A \emph{Fuchsian group of hyperbolic space} is a discrete group of orientation-preserving isometries leaving globally invariant a
totally geodesic surface, denoted by $P_{\mathbb{H}^2}$, on which it acts cocompactly and without fixed points. 
In \cite{artrealisationhyp} it is proved that:
\begin{theorem}\label{realisationhyp}
A hyperbolic metric with conical singularities of positive singular curvature on a compact surface $S$ of genus $>1$ is realised by a unique (up to global isometries) convex Fuchsian polyhedron in hyperbolic space.
\end{theorem}
By the Gauss--Bonnet formula, we know that there doesn't exist other constant curvature metrics with conical singularities of constant sign singular curvature on surfaces of genus $>1$ than the ones described in  Theorem~\ref{realisation lorentz} and Theorem~\ref{realisationhyp}.

A \emph{parabolic polyhedron} is a polyhedral surface of hyperbolic or de Sitter space invariant under the action of a group of isometries acting freely cocompactly  on a horosphere. Convex parabolic polyhedra provide constant curvature metrics with conical singularities with constant sign singular curvature on the torus. We think that every such metric is realised by a unique convex parabolic polyhedron. The hyperbolic case is done in \cite{metrictorus}. For all genus the hyperbolic results were extended to metrics with cusps and complete ends of infinite area in \cite{CompHyp}.

If the parabolic de Sitter case is done, that would answer the following question. Let $S$ be a compact surface with fundamental group $\Gamma$.
\begin{gstatement}
Can each constant curvature $K$ metric with conical singularities with constant sign singular curvature $\epsilon \in \{-,+\}$ on  $S$ be realised in $M_K^{\epsilon}$ by a unique (up to global isometries) convex  polyhedral surface invariant under the action of a representation of $\Gamma$ in a group of isometries acting freely cocompactly on a totally umbilic surface? (with a condition 
on the lengths of contractible geodesics in the cases $K=1,\epsilon=-$).
\end{gstatement}

\subsection{Andreev's Theorem.} 
Theorem~\ref{thmrivhod} can be seen as a generalisation of the famous Andreev's Theorem about compact hyperbolic polyhedra with acute dihedral angles \cite{andreev2,andreevcorrected}. It is proved in  \cite{hodgsonandreev} and it seems that the genus does not intervene in this proof. It follows that  the part 1) of Theorem~\ref{realisation lorentz} would be seen as a generalisation of the Andreev's Theorem for surfaces of genus $>1$.

\subsection{Dual statements.}

Using the duality between hyperbolic space and de Sitter space (see Subsection~\ref{subsec:duality}), part 1) of Theorem~\ref{realisation lorentz} 
can be reformulated as a purely hyperbolic statement:
\begin{Alphatheorem}\label{enonce dual du thm realisation hyperbolique}
Let $S$ be a compact surface of genus $>1$ with a spherical metric $h$ 
with conical singularities of negative singular  curvature 
such that 
its closed contractible geodesics have lengths   $>2\pi$.

There exists a unique (up to global isometries) convex Fuchsian polyhedron in the hyperbolic space such that   its dual metric is  isometric to $h$ (up to a quotient).
\end{Alphatheorem}

We could make the same statement in the anti-de Sitter case, which is its own dual.

\subsection{Hyperbolic manifolds with polyhedral boundary.}

Take a convex Fuchsian polyhedron $(P,F)$ of the hyperbolic space and consider the Fuchsian polyhedron
$(P',F)$ obtained by a symmetry relative to the plane $P_{\mathbb{H}^2}$ (the one fixed by the Fuchsian group action). Next cut the hyperbolic space
along $P$ and $P'$, and keep the  component bounded by $P$ and $P'$. The quotient of this manifold by the acting group $F$ is a kind of
hyperbolic manifold called \emph{Fuchsian manifold} (with
convex polyhedral boundary): they are compact hyperbolic manifolds with convex boundary with an
isometric involution fixing a hyperbolic surface (the symmetry relative to
$P_{\mathbb{H}^2}/F$). All the Fuchsian manifolds can be
obtained in this way: the lifting to the universal covers of the canonical embedding of a component of the boundary in the Fuchsian manifold 
gives a Fuchsian polyhedron of the hyperbolic space.
Theorem~\ref{enonce dual du thm realisation hyperbolique} says
exactly that for a choice of a (certain kind of) spherical  metric $g$ on  the boundary, there
exists a unique hyperbolic metric on the Fuchsian manifold such that the
dual metric of the induced metric of the boundary is isometric to $g$:
\begin{Alphatheorem}
The hyperbolic metric on a Fuchsian manifold with convex polyhedral boundary is uniquely determined by the dual  metric of its boundary.
\end{Alphatheorem}

This can be generalised as the following question.
Let $M$ be a compact connected manifold with boundary $\partial M$ of dimension 3, which admits a complete hyperbolic convex cocompact metric. We know that the dual metric of the induced metric on $\partial M$ is a spherical metric with conical singularities of negative singular curvature such that the lengths of 
its closed contractible geodesics are $>2\pi$.
\begin{gstatement}\label{conjecture variete2}
Is each such dual metric   on $\partial M$  induced on $\partial M$ by a unique (up to isometries) hyperbolic metric on $M$?
\end{gstatement}

Theorem~\ref{thmrivhod} is a positive answer to this question in the case where $M$ is the ball.
The analogous of Question~\ref{conjecture variete2} for the induced metric on the boundary is stated in \cite{artrealisationhyp}, and the Fuchsian particular case corresponds to Theorem~\ref{realisationhyp}.
Both  questions in the case where the boundary is smooth and strictly convex have been positively answered in \cite{Schconvex}.

It is also possible to do analogous statements for  ``anti-de Sitter manifolds with convex boundary'', related to  so-called \emph{maximal globally hyperbolic anti-de Sitter manifolds}, see \emph{e.g.} \cite{mess,Schgren}. The part 3) of Theorem~\ref{realisation lorentz} would describe a particular case of this statement.

\subsection{Convention.} \label{subsec:convention}
In all the paper, we call \emph{length} of a time-like geodesic the imaginary part of the ``distance'' between the endpoints of the geodesic. We will also call \emph{distance} between two points joined by a time-like geodesic the length (in the sense we have just defined) of the geodesic. It follows that distances and lengths will be real numbers. In this mind, for a time-like vector $X$, we denote the modulus of its ``norm'' by $\norm{X}$ instead of $\vert \norm{X} \vert$.

\subsection{Acknowledgements.}
The material in this paper, together with \cite{artrealisationhyp}, is a part of my PhD thesis  \cite{these} under the direction of
Bruno Colbois and Jean-Marc Schlenker. For that reason, they played a
crucial part in the working out of these results. 
I also want to thank Christophe Bavard, Michel Boileau and Marc Troyanov for their useful comments.
The paper was also crucially improved and simplified thanks to remarks from an anonymous referee.

This paper was first announced under the title ``\emph{Polyhedral realisation of metrics with conical singularities on compact surfaces in Lorentzian space-forms}''.

The author was first supported by the Mathematic Institute of the University of Neuch\^atel and the Laboratoire \'Emile Picard of the University Paul Sabatier (Toulouse) and was then partially supported by the Schweizerischer Nationalfonds 200020-113199/1. 

%%%%%%%%%%%%%%%%%%%%%%%%%%%%%%%%
%%%%%%%%%%%%%%%%%%%%%%%%%%%%%%%%

\section{Backgrounds}\label{section backgrounds}

\subsection{Projective models.}\label{modelesprojectifs}

For this subsection and next, details can be found in \cite{RivHod,schhilb,schlorentz}. 
 The notation  $\mathbb{R}_p^n$ means $\mathbb{R}^n$ endowed with a non-degenerate bilinear form of signature $(n-p,p)$ --- read $(+,-)$. In particular we denoted $\mathbb{R}_1^3$ as $M_0^-$. The ``norm'' associated to $\mathbb{R}_p^n$ is denoted by $\norm{.}_p$. Recall that Riemannian and Lorentzian space-forms can be seen as pseudo-spheres in flat spaces:
\begin{eqnarray*}
\ &&\mathbb{H}^3=\{x\in\mathbb{R}^4_1 \vert \norm{x}_1^2=-1, x_4>0\}, \\
\ &&M_1^-=\{x\in\mathbb{R}^4_1 \vert \norm{x}_1^2=1\}, \\
\ &&M^-_{-1}=\{x\in\mathbb{R}^4_2 \vert \norm{x}_2^2=-1\}
\end{eqnarray*}
(note that the negative directions are always the last ones).

In each $M_K^-$ we choose the following particular point $c_K$:
\begin{eqnarray*} 
\ &&c_1=(1,0,0,0)\in \mathbb{R}^4_1, \\
\ &&c_0=(0,0,0)\in M^-_0,\\
\ &&c_{-1}=(0,0,0,1)\in\mathbb{R}^4_2.
\end{eqnarray*}

With this definition of $c_1$, for example, the hyperbolic surface $O_1$ is given by the intersection in $\mathbb{R}^4_1$ of the pseudo-sphere $M^-_1$ and the hyperplane $\{x_1=\sqrt{2}\}$.
We will denote by $\Omega_K$ the future-cone of the point $c_K$ in $M_K^-$ for $K=0,1$. For $K=-1$, we denote by  $\Omega_{-1}$ the intersection of the future-cone of $c_{-1}$ with $\{x_4\geq 0\}$.

\begin{figure}[ht] \begin{center}
\begin{picture}(0,0)%
\includegraphics{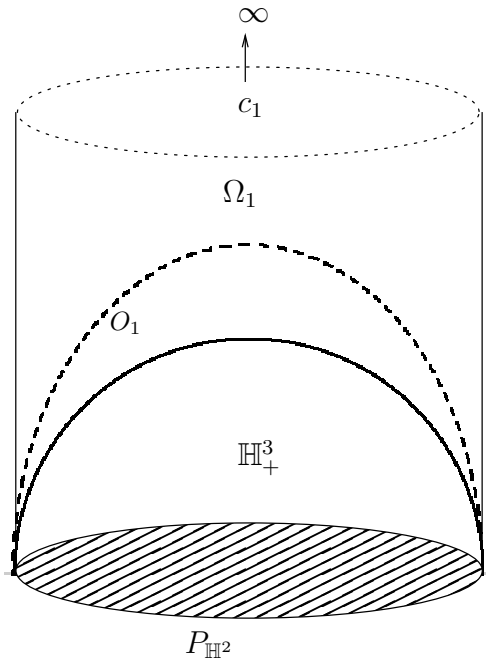}%
\end{picture}%
\setlength{\unitlength}{4144sp}%
\begingroup\makeatletter\ifx\SetFigFont\undefined%
\gdef\SetFigFont#1#2#3#4#5{%
  \reset@font\fontsize{#1}{#2pt}%
  \fontfamily{#3}\fontseries{#4}\fontshape{#5}%
  \selectfont}%
\fi\endgroup%
\begin{picture}(2901,4024)(11161,-4229)
\put(12241,-4156){\makebox(0,0)[lb]{\smash{{\SetFigFont{12}{14.4}{\rmdefault}{\mddefault}{\updefault}{\color[rgb]{0,0,0}$P_{\mathbb{H}^ 2}$}%
}}}}
\put(11791,-2221){\makebox(0,0)[lb]{\smash{{\SetFigFont{10}{12.0}{\rmdefault}{\mddefault}{\updefault}{\color[rgb]{0,0,0}$O_1$}%
}}}}
\put(12556,-916){\makebox(0,0)[lb]{\smash{{\SetFigFont{12}{14.4}{\rmdefault}{\mddefault}{\updefault}{\color[rgb]{0,0,0}$c_1$}%
}}}}
\put(12556,-376){\makebox(0,0)[lb]{\smash{{\SetFigFont{12}{14.4}{\rmdefault}{\mddefault}{\updefault}{\color[rgb]{0,0,0}$\infty$}%
}}}}
\put(12466,-1456){\makebox(0,0)[lb]{\smash{{\SetFigFont{12}{14.4}{\rmdefault}{\mddefault}{\updefault}{\color[rgb]{0,0,0}$\Omega_1$}%
}}}}
\put(12556,-3031){\makebox(0,0)[lb]{\smash{{\SetFigFont{12}{14.4}{\rmdefault}{\mddefault}{\updefault}{\color[rgb]{0,0,0}$\mathbb{H}^3_+$}%
}}}}
\end{picture}%

\caption{Klein projective model of hyperbolic and de Sitter spaces\label{spheredesitter}.}
\end{center}
\end{figure}

The \emph{Klein projective models} of the hyperbolic and  de Sitter spaces are the images of both spaces under a projective map given by the projection from the origin of $\mathbb{R}^4_1$ onto the hyperplane $\{x_4=1\}$: $x\mapsto x/x_4$ (the hyperplane $\{x_4=1\}$ is naturally endowed with a Euclidean structure in $\mathbb{R}^4_1$).
The hyperbolic space is sent to the open unit ball and the de Sitter space is sent to the exterior of the closed unit ball (actually, to have a diffeomorphism, we involve in this projection only the upper part of the de Sitter space given by $\{x_4\geq 0\}$, but this is not restrictive anyway because all the surfaces we will consider will be contained in $\Omega_1$, itself contained in the upper part of the de Sitter space). 
The projective map sends the point $c_1$ to infinity, and $\Omega_1$ is sent to an infinite cylinder with basis a unit disc centered at the origin, and this disc corresponds to the hyperbolic plane dual to $c_1$ (see below), denoted by $P_{\mathbb{H}^2}$. We denote by $\mathbb{H}^3_+$ the intersection of the hyperbolic space with the half-plane defined by $P_{\mathbb{H}^2}$ and containing $c_1$.  In this model $O_1$ (as well as all surfaces in $\Omega_1$ at constant distance from $c_1$) is a half-ellipsoid with same boundary as   $P_{\mathbb{H}^2}$, \emph{i.e.} the  horizontal unit circle, see Figure~\ref{spheredesitter}. The \emph{boundary at infinity} of both hyperbolic and de Sitter space in this model is the unit sphere.

\begin{figure}[ht] \begin{center}
\begin{picture}(0,0)%
\includegraphics{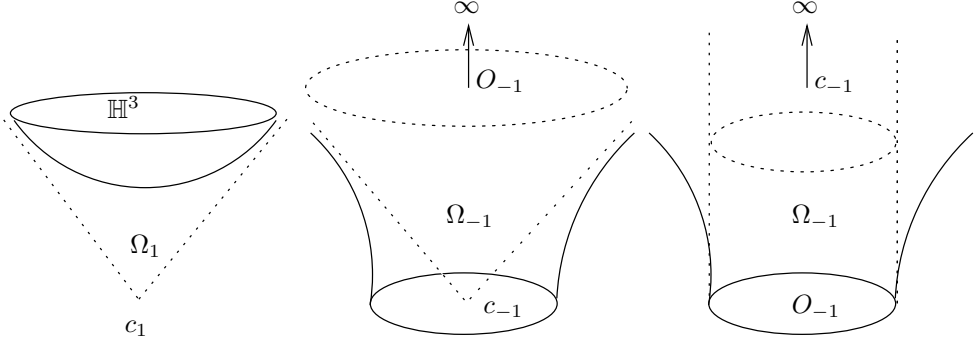}%
\end{picture}%
\setlength{\unitlength}{4144sp}%
\begingroup\makeatletter\ifx\SetFigFont\undefined%
\gdef\SetFigFont#1#2#3#4#5{%
  \reset@font\fontsize{#1}{#2pt}%
  \fontfamily{#3}\fontseries{#4}\fontshape{#5}%
  \selectfont}%
\fi\endgroup%
\begin{picture}(5827,2106)(304,-3607)
\put(3196,-3436){\makebox(0,0)[lb]{\smash{{\SetFigFont{10}{12.0}{\rmdefault}{\mddefault}{\updefault}{\color[rgb]{0,0,0}$c_{-1}$}%
}}}}
\put(3016,-1636){\makebox(0,0)[lb]{\smash{{\SetFigFont{10}{12.0}{\rmdefault}{\mddefault}{\updefault}{\color[rgb]{0,0,0}$\infty$}%
}}}}
\put(2971,-2896){\makebox(0,0)[lb]{\smash{{\SetFigFont{10}{12.0}{\rmdefault}{\mddefault}{\updefault}{\color[rgb]{0,0,0}$\Omega_{-1}$}%
}}}}
\put(3151,-2086){\makebox(0,0)[lb]{\smash{{\SetFigFont{10}{12.0}{\rmdefault}{\mddefault}{\updefault}{\color[rgb]{0,0,0}$O_{-1}$}%
}}}}
\put(5041,-2896){\makebox(0,0)[lb]{\smash{{\SetFigFont{10}{12.0}{\rmdefault}{\mddefault}{\updefault}{\color[rgb]{0,0,0}$\Omega_{-1}$}%
}}}}
\put(5041,-1636){\makebox(0,0)[lb]{\smash{{\SetFigFont{10}{12.0}{\rmdefault}{\mddefault}{\updefault}{\color[rgb]{0,0,0}$\infty$}%
}}}}
\put(5176,-2086){\makebox(0,0)[lb]{\smash{{\SetFigFont{10}{12.0}{\rmdefault}{\mddefault}{\updefault}{\color[rgb]{0,0,0}$c_{-1}$}%
}}}}
\put(5041,-3436){\makebox(0,0)[lb]{\smash{{\SetFigFont{10}{12.0}{\rmdefault}{\mddefault}{\updefault}{\color[rgb]{0,0,0}$O_{-1}$}%
}}}}
\put(1050,-3547){\makebox(0,0)[lb]{\smash{{\SetFigFont{10}{12.0}{\rmdefault}{\mddefault}{\updefault}{\color[rgb]{0,0,0}$c_1$}%
}}}}
\put(946,-2266){\makebox(0,0)[lb]{\smash{{\SetFigFont{10}{12.0}{\rmdefault}{\mddefault}{\updefault}{\color[rgb]{0,0,0}$\mathbb{H}^3$}%
}}}}
\put(1081,-3076){\makebox(0,0)[lb]{\smash{{\SetFigFont{10}{12.0}{\rmdefault}{\mddefault}{\updefault}{\color[rgb]{0,0,0}$\Omega_1$}%
}}}}
\end{picture}%
\caption{Minkowski projective models for hyperbolic and de Sitter space (left) and for anti-de Sitter space (center). On the right is the second Minkowski projective model for anti-de Sitter space \label{modeleprojectifMink}}
\end{center}
\end{figure}

Using the projection onto the hyperplane $\{x_1=1\}$ (naturally isometric to $M^-_0$), we can define another projective model for both hyperbolic and de Sitter spaces, called  \emph{Minkowski projective model}. In this model,  $\mathbb{H}^3_+$ is sent to the interior of the (unitary) upper branch of the two-sheeted hyperboloid and $P_{\mathbb{H}^2}$ is sent to infinity. The image of (a part of) the de Sitter space lies between the light-cone and this upper branch. The point $c_1$ is sent to the origin and $\Omega_1$ is sent to the complementary in $ \Omega_0$ of the interior of the upper branch of the (unitary) two-sheeted hyperboloid, see Figure~\ref{modeleprojectifMink} left. We denote by $\varphi_1$ this map.  The \emph{boundary at infinity} of both hyperbolic and de Sitter space in this model is  the upper branch of the two-sheeted hyperboloid.

There exists also a Minkowski projective model for the anti-de Sitter space, using the projection onto the hyperplane $\{x_4=1\}$. The anti-de Sitter space is sent to the interior of the (unitary) one-sheeted hyperboloid, which is the \emph{boundary at infinity} of anti-de Sitter space in this model. The point
$c_{-1}$ is sent to the origin, $\Omega_{-1}$ is sent to $\Omega_0$ and $O_{-1}$ is sent to infinity, see Figure \ref{modeleprojectifMink} center. 
We denote by $\varphi_{-1}$  this map.

For the anti-de Sitter space, there exists another projection onto Minkowski space, that we will call \emph{second Minkowski projective model}. It suffices to project onto $\{x_3=1\}$ instead of $\{x_4=1\}$. Then  $c_{-1}$ is sent to infinity, $O_{-1}$ is sent to the horizontal disc, and  surfaces at constant distance from $c_{-1}$ are sent to half-ellipsoids with boundary  the unit circle in the horizontal plane and  $\Omega_{-1}$ corresponds to the interior of the half-cylinder above $O_{-1}$, see Figure \ref{modeleprojectifMink} right. The  future-cone of $c_{-1}$ is the interior of the whole cylinder.  The boundary at infinity is still the one-sheeted hyperboloid.

All these projective maps send  geodesics onto geodesics, and hence convex sets onto convex sets. 
In the projective models of de Sitter space, the geodesics  are space-like if they don't meet the boundary at infinity, light-like if they are tangent to the boundary at infinity and time-like if they meet the boundary at infinity. 
 In the projective models of anti-de Sitter space, the geodesics  are space-like if they  meet the boundary at infinity, light-like if they are tangent to the boundary at infinity and time-like if they don't meet the boundary at infinity.
Hence it is easy to draw convex space-like polyhedral surfaces in these models.

\subsection{Description of Fuchsian polyhedra.}

By definition, the vertices of a convex Fuchsian polyhedron are the union of finitely many orbits of points lying on surfaces at constant distance from $c_K$, and, for $K=0,1$, these surfaces 
lie inside $\Omega_K$. Moreover, for $K=0,1$, $c_K$ is in the concave side of such surfaces, then $c_K$ is also in the concave part of the convex Fuchsian polyhedra.

\begin{lemma}\label{region des pol fuch ads}
Up to an isometry, a convex Fuchsian polyhedron in the anti-de Sitter space lies between $c_{-1}$ and 
$O_{-1}$, \emph{i.e.} it is contained in $ \Omega_{-1}$, and $c_{-1}$ is in its concave part.
\end{lemma}
\begin{proof}
In the second Minkowski projective model, the vertices of the convex Fuchsian polyhedron are lying on half-ellipsoids having same boundary as the horizontal unit circle, which represents $O_{-1}$. By convexity, all these half-ellipsoids must be on the same side of $O_{-1}$, and up to an isometry, the side is chosen to be the one containing $c_{-1}$.
\end{proof}

As $O_K$ is at constant distance from $c_K$, the time-like geodesics from $c_K$ are orthogonal to $O_K$, then 
they define an orthogonal projection of the future-cone of $c_K$ onto $O_K$ that will be denoted by $p_K$.  Moreover
\begin{lemma}
The maps $p_K$ are homeomorphisms between each convex Fuchsian polyhedron $P$ and $O_K$.
\end{lemma}
\begin{proof}
We will prove it in the de Sitter case, the others cases follow immediately using the projective maps. In  the Klein projective model, $p_1$ corresponds to the vertical projection, and 
then the orthogonal projection of $P$ onto the horizontal plane is one-to-one, as  the convex hull of $P$ is the union of $P$ and of the closed unit disc of the horizontal plane. It follows that this horizontal projection is a homeomorphism.  
Moreover, as $O_1$ is a  convex cap,  its orthogonal projection onto the horizontal plane
is a homeomorphisms. And it suffices to compose these projections to get $p_1$.
\end{proof}

\begin{lemma}\label{fuchsien dans conelum est espace}
A convex polyhedral surface  contained in $\Omega_0$ or $\Omega_{-1}$  and invariant under the action of a Fuchsian group is a convex Fuchsian polyhedron, \emph{i.e.} it has space-like faces and space-like edges. 
\end{lemma}
 \begin{proof}
We will prove that the faces are space-like, that will imply that the edges are space-like, as they lie in space-like planes.  
Suppose that the polyhedral surface has a non-space-like face $H$.
We denote by $S(x)$ the surface at constant distance from $c_K$ which contains a vertex $x$ belonging to  $H$. 
The images of $x$ by the action of the Fuchsian group $F$ are all lying on $S(x)$. Moreover, up to multiply the metric by a constant, the induced metric on $S(x)$ is hyperbolic, and by cocompactness the limit set of $F$  corresponds to the entire boundary at infinity of $S(x)$.

For the anti-de Sitter case, it is easy to see (in particular in the second Klein projective model) that a non-space-like plane meeting $S(x)$ separates its boundary at infinity in two components, each one containing an infinite number of points. By convexity, all the images of $x$ under the  action of $F$ must stay in the same side of $H$, that is impossible as all the points of the boundary at infinity must be reached. 

Let consider now the Minkowski case. If $H$ belongs to a time-like affine hyperplane, we can use the same argument as above. 
Things are more subtle if $H$ belongs to a light-like affine hyperplane. In this case its intersection with $S(x)$ is a horocycle, it follows that all the images of $x$ by the action of $F$ must be ``outside'' the horocycle (\emph{i.e.} in its non-convex part in the Klein projective model). To get a contradiction with the convexity, it is enough to show that if $(x_k)_k$ is a sequence such that $x_k:=f_kx$, $f_k\in F$, converging to the ``center'' of the horocycle on the boundary at infinity, then the sequence must meet the ``inner part'' of the horocycle (\emph{i.e.} its convex part in the Klein projective model). But this is a well-known fact of the  geometry of the hyperbolic plane, as $F$ is a cocompact group containing only hyperbolic isometries. 
\end{proof} 

Such a property is false in the de Sitter case. The fact is that $\varphi_1$ is not a bijection from $\Omega_1$ to $\Omega_0$. Or equivalently, if $G$ is a geodesic plane and $S$ is a surface of $\Omega_{1}$ at constant distance from $c_1$, the intersection of $S$ and $G$ can be a circle (\emph{i.e.} $G$ does not meet the boundary at infinity of $S$) and $G$ can be not space-like, that is impossible in Minkowski and anti-de Sitter cases.

\begin{lemma}\label{proj envoie fuchs sur fuchs}
The projective map $\varphi_{-1}$ (resp. $\varphi_1$)
sends convex Fuchsian polyhedra  to convex Fuchsian polyhedra.
\end{lemma}
\begin{proof}
We already know that $\varphi_{-1}$ (resp. $\varphi_1$) sends  convex polyhedral surfaces contained in $\Omega_{-1}$ (resp. $\Omega_1$) to  convex polyhedral surfaces of $\Omega_0$. We will see that they also
act on the representations of $\Gamma$. This, together with Lemma~\ref{fuchsien dans conelum est espace}, will prove the statement.

We denote by $\mathrm{\mathrm{Isom}}_{c_1}(M_1^-)$ the subgroup of $\mathrm{Isom}^+_+(M_1^-)$ which fix the point  $c_1$. The fact is that: 
the projective map $\varphi:M_1^-\longrightarrow M^-_0$
induces an isomorphism $G:\mathrm{Isom}_{c_1}(M_1^-) \longrightarrow \mathrm{Isom}_{c_0}(M^-_0)$ which commutes with 
the projective map, that is, if $f\in \mathrm{Isom}_{c_{1}}(M_1^-)$, then
\begin{equation}
\varphi(f(x))=G(f)(\varphi(x)).\nonumber
\end{equation}
The projective map from $M_{-1}^-$ to $M^-_0$  has the same property. 
Now we prove this fact.

In the  definition of $\varphi_1$, the Minkowski space $M_0^-$ 
is seen as the intersection of $\mathbb{R}^{4}_1$ with the hyperplane  $\{x_1=1\}$. 
 It allows us to extend the isometries of the Minkowski space of dimension $3$: an isometry of $M_0^-$
sending $x$ to $y$ can be extended to an isometry of $\mathbb{R}^{4}_1$ 
which sends the point  $(t,tx)$ to the point  $(t,ty)$. Furthermore, these isometries of  $\mathbb{R}^{4}_1$ preserve 
the de Sitter space: there restrictions to the hyperboloid are isometries of the de Sitter space. It follows that the isometries of 
$M_0^-$ which fix the origin correspond to isometries of de Sitter space which fix the point $c_1$.

It is the same thing in the anti-de Sitter case, isometries which send $x$ to $y$ are extended to isometries 
sending $(tx,t)$ to $(ty,t)$. 
The properties of commutations are then obvious. And by construction, a cocompact group is sent to a cocompact group.
\end{proof}

\begin{lemma}\label{lem:outlightcone}
  Each convex Fuchsian polyhedron of $M_K^-$, $K\in\{-1,0,1\}$, stays out of the light-cone of its vertices.
\end{lemma}
\begin{proof}
 Let $(P,F)$ be a convex Fuchsian polyhedron and $x$ a vertex of $P$. As $c_K$ belongs to the concave part of $P$, it is clear that the part of the light-cone of $x$ such that its interior contains $c_K$ never meets $P$. We denote by $C(x)$ the other part of the light-cone, namely the one bounding the future-cone of $x$. The orbit of $x$ under the action of $F$ lies on a surface $S(x)$ at constant distance from $c_K$. Suppose that a point $y$ of $P$ meets $C(x)$. The point $y$ lies on a space-like plane $H$, $H$ containing a face of $P$. In de Sitter and Minkowski cases, it is not hard to see that  the intersection of $H$ with $S(x)$ is a circle, splitting $S(x)$ into a compact and a non-compact part. In anti-de Sitter case, the intersection of $S(x)$ with a space-like plane can be different from a circle, but as here the plane contains a face of a Fuchsian polyhedron,  the intersection must be a circle (the argument is the same as the proof of Lemma~\ref{fuchsien dans conelum est espace} for the Minkowski case, depending on how the plane meets the boundary at infinity of $S(x)$). Moreover  
 the light-cone of $y$  has a common direction with $C(x)$,  and this implies that $x$ lies in the compact part. It follows by convexity that the orbit of $x$  remains in this compact part, that is impossible. 
\end{proof}

\subsection{Duality.}\label{subsec:duality}

There exists a well-known duality between hyperbolic space and de Sitter space  \cite{coxeter,Thurcour1,rivinthese,RivHod,Thurlivr1,schhilb}. It corresponds in the Klein projective model to the classical projective duality with respect to the light-cone (\emph{i.e.} the unit sphere in the Euclidean space). In $\mathbb{R}^4_1$, the dual of a unit vector (a point on a pseudosphere) is its orthogonal hyperplane and the dual of a plane is its orthogonal plane. It follows that the dual of a face is a vertex, the dual of a vertex is a face and the dual of an edge is an edge, the dual of a convex polyhedral surface is a convex polyhedral surface, and the duality is an involution. 

The dual of a convex polyhedral surface of the hyperbolic space is a space-like convex polyhedral surface of the de Sitter space.  Polyhedral surface of the de Sitter space obtained by this way are called \emph{polyhedral surfaces of hyperbolic type}. There exists (space-like convex) polyhedral surfaces in the de Sitter space which are not of hyperbolic type (their dual is not contained in the hyperbolic space), see \cite{schlorentz}.

\begin{lemma}
A convex Fuchsian polyhedron of the de Sitter space is of hyperbolic type.
\end{lemma}
\begin{proof}
Using a projective description of the duality in term of cross-ratio, it is easy to check that the dual of a surface in $\Omega_1$ at constant distance from $c_1$ is a hyperbolic surface  at constant distance from the hyperbolic plane dual to $c_1 $ \cite{Schnoncpt}. As a convex Fuchsian polyhedron lies between two surfaces at constant distance from $c_1$, it follows that its dual lies between two hyperbolic surfaces, then it lies entirely in the hyperbolic space.
\end{proof}

And the lengths of the closed geodesics for the induced metric on a convex polyhedral surface of hyperbolic type are $>2\pi$ \cite{RivHod}, that explains the additional condition in the part 1) of Theorem \ref{realisation lorentz}.

Note that the proof of this Lemma also says that the dual of a convex Fuchsian polyhedron of the hyperbolic space is contained inside the future-cone of $c_1$.

\begin{lemma}\label{dual fuchsien est dual}
The dual of a convex Fuchsian polyhedron is a convex Fuchsian polyhedron.
\end{lemma}
\begin{proof}
A Fuchsian polyhedron in the hyperbolic space is invariant under the action of a group $F$ of isometries which leaves invariant  $P_{\mathbb{H}^2}$, and this one is given by the intersection of the hyperbolic space with a time-like hyperplane $V$ in $\mathbb{R}^4_1$. The group $F$ is  given by orientation preserving and time orientation preserving isometries of Minkowski space leaving invariant $V$. These isometries also fix the unit space-like  vector $v$ normal to $V$ (which corresponds to the point $c_1$ of the de Sitter space), and they also fix the de Sitter space. Moreover these isometries act cocompactly on all the hyperplanes orthogonal to $V$, in particular the one which defines $O_1$. It follows that the restrictions of these isometries to the de Sitter space are Fuchsian isometries. The converse holds in the same manner.
\end{proof}

Moreover,  the dual metric of a convex polyhedral surface in hyperbolic space is isometric to the metric induced on its dual \cite{RivHod,CharDav1}, that explains Theorem~\ref{enonce dual du thm realisation hyperbolique}.

%%%%%%%%%%%%%%%%%%%%%%%
\subsection{Teichm\"uller space}\label{Teichspace}

\subsubsection{Z-V-C  coordinates for Teichm\"uller space.}

For more details about  Z-V-C coordinates (Z-V-C  stands for Zieschang--Vogt--Coldewey, \cite{ZVC}) we refer to \cite[6.7]{Buser}.
\begin{definition}
Let $g\geq 2$. A (geodesically convex) polygon of the hyperbolic plane with edges  (in the direct order)  $b_1,b_2,\overline{b}_1,\overline{b}_2,b_3,b_4,\ldots,\overline{b}_{2g}$ and with interior angles $\theta_1,\overline{\theta}_1,\ldots,\theta_{2g},\overline{\theta}_{2g}$ is called   \emph{(normal) canonical} if, with $l(c)$ the length of the geodesic $c$,
\begin{itemize}
\item[i)] $l(b_k)=l(\overline{b}_k)$, $\forall k$;
\item[ii)] $\theta_1+\ldots+\overline{\theta}_{2g}=2\pi$;
\item[iii)] $\theta_1+\theta_2=\overline{\theta}_1+\theta_2=\pi$.
\end{itemize}
Two canonical polygons $P$ and $P'$ with edges $b_1,\ldots,\overline{b}_{2g}$ and  $b'_1,\ldots,\overline{b}'_{2g}$ are said
\emph{equivalent} if there exists an isometry from $P$ to $P'$ such that the edge  $b_1$ is sent to the edge $b'_1$ and $b_2$ is sent to $b'_2$.
\end{definition}

If we identify the edges $b_i$ with the edges $\overline{b}_i$, we
get a compact  hyperbolic surface of genus $g$. This surface could
also be written $\mathbb{H}^2 / F$, where $F$ is the sub-group of
$\mathrm{PSL}(2,\mathbb{R})=\mathrm{Isom}^+(\mathbb{H}^2)$ generated by the translations along the edges
$b_i$ (the translation length is the length of $b_i$).
The interior of the polygon is a fundamental domain for the action
of $F$. This leads to a description of the Teichm\"uller space
$\mathrm{T}_g$:
\begin{proposition}[{\cite[6.7.7]{Buser}}]
Let $\mathrm{P}_g$  be the set of equivalence classes of canonical polygons. An element of  $\mathrm{P}_g$
is described by the $(6g-6)$ real numbers (the Z-V-C coordinates):
\begin{equation}
(b_3,\ldots,b_{2g},\theta_3,\overline{\theta}_3,\ldots,\theta_{2g}, \overline{\theta}_{2g}).\nonumber
\end{equation}
Endowed with this topology, $\mathrm{P}_g$ is in analytic bijection with $\mathrm{T}_g$.
\end{proposition}

\subsubsection{Fenchel--Nielsen coordinates for Teichm\"uller space.}

A compact hyperbolic surface $S$ of genus $g$ can be described as a gluing of pants. 
Such a gluing leads to the choice of $(6g-6)$ real numbers: the lengths of the geodesics along which we glue the pants and the angles of the twists. These numbers (the Fenchel--Nielsen coordinates) describe the  Teichm\"uller space of $S$ \cite{Buser}. 

We denote by $l(x)$ the length of the curve $x$. It is possible to compute the twist parameters knowing the lengths of certain geodesics:
\begin{proposition}[{\cite[3.3.12]{Buser}}]\label{resultat buser}
Let  $\gamma$ be the geodesic along which are glued two pants. We denote by  $\delta$ the geodesic along which we can cut the $X$-piece to get the other pant decomposition than the one given by $\gamma$. 

We do a twist with parameter $\alpha$ around  $\gamma$, and we denote by  $\delta^{\alpha}$ the geodesic which is in the homotopy class of the image of $\delta$ by the twist. Then there exists real analytic functions of  $l(\gamma)$ $u$ and $v$ ($v>0$)   such that
\begin{equation}\label{equation FN long geod orth}
\cosh (\frac{1}{2} l(\delta^{\alpha})) = u+ v \cosh (\alpha l(\gamma)).
\end{equation}
\end{proposition}

Let $(S_k,h_k)_k$ be a sequence of (equivalence classes) of hyperbolic surfaces in the Teichm\"uller space of $S$. The metric on $S$ is  $h_0$. We denote by $f_k$ a homeomorphism between $S$ and $S_k$,  and by $\mathrm{L}_{h_k}(\gamma)$ the length of the geodesic corresponding to the element $(f_k)_*(\gamma)$ of the fundamental  group of $S_k$ ($\gamma\in \pi_1(S)$) for the metric  $h_k$.

\begin{lemma}\label{lemme teichmuller minor geod}
If there exists a constant  $c>0$ such that, for all $\gamma\in\pi_1(S)$, $\mathrm{L}_{h_k}(\gamma)\geq \frac{1}{c} \mathrm{L}_{h_0}(\gamma)$, then $(h_k)_k$ converges (up to extract a subsequence).
\end{lemma}
\begin{proof}
We will prove  that in this case  $(h_k)_k$ is contained inside a compact of the   Teichm\"uller space using the Fenchel--Nielsen coordinates. 
First, the lengths of the geodesics along which we cut up to obtain the pants are bounded from above: if not, because of the  Gauss--Bonnet Formula, it would exist another geodesic with decreasing length   (otherwise the area of the surface may become arbitrary large), that is impossible as the lengths are bounded from below.

It remains to check that the twist parameters are bounded. Take a geodesic $\gamma$ along which a twist is done. 
We suppose now that the twist parameter $\alpha_k$ associated to $\gamma$ becomes arbitrary large  (up to consider a $k$ sufficiently large and up to extract a sub-sequence, in regard of the considerations above we can suppose that the length of  $\gamma$ doesn't change). We look at the  geodesic $\delta$ defined as in Proposition \ref{resultat buser} for the metric $h_k$. If now we do the twists in the counter order to go from  $h_k$ to $h_0$ (\emph{i.e.} we do twists with parameters $-\alpha_k$), then by (\ref{equation FN long geod orth}), $\mathrm{L}_{h_0}(\delta)$ is arbitrariness  larger than $\mathrm{L}_{h_k}(\delta)$. 
\end{proof}

\subsubsection{Geodesic lengths coordinates for Teichm\"uller space.}

With the same notations than above, here is a result that can be found \emph{e.g.} as \cite[Expos\'e 7, Proposition 5]{travauxthurstonsurfaces}:
\begin{lemma}\label{lemme teichmuller major geod}
There exists a finite number of elements $(\gamma_1,\ldots,\gamma_n)$ of the fundamental group of $S$ such that, if for all $i$  $\mathrm{L}_{h_k}(\gamma_i)$ is uniformly bounded from above, then $(h_k)_k$ converges (up to extract a subsequence).
\end{lemma}

%%%%%%%%%%%%%%%%%%%%%%%%%%%%%%%%%%%%%%%%%%%%%%%%%%%%%%%%%%%%%%%%%%%%
%%%%%%%%%%%%%%%%%%%%%%%%%%%%%%%%%%%%%%%%%%%%%%%%%%%%%%%%%%%%%%%%%%%%%
\section{Fuchsian infinitesimal rigidity}\label{rigidite lorentz}

\subsection{Fuchsian polyhedral embeddings.}

To define the ``Fuchsian infinitesimal rigidity'' we  need to   describe Fuchsian polyhedra as polyhedral embeddings. 

\begin{definition}
A (space-like) \emph{polyhedral embedding} of a surface $S$ into $M_K^-$ is  a cellulation of $S$ together with a homeomorphism from 
$S$ to a (space-like) polyhedral surface of  $M_K^-$, sending polygons of the cellulation
to (space-like) geodesic polygons of  $M_K^-$. 

A \emph{Fuchsian polyhedral  embedding}
 in   $M_K^-$ is a triple
 $(S,\phi,\rho)$, where
\begin{itemize}
\item $S$ is a compact surface  of genus $>1$,
\item $\phi$ is a polyhedral embedding
of the universal cover $\widetilde{S}$ of $S$ into $M_K^-$,
\item  $\rho$ is a representation of the fundamental group $\Gamma$ of  $S$ into $\mathrm{Isom}_+^+(M_K^-)$,
\end{itemize}
such that $\phi$ is \emph{equivariant} under the action of $\Gamma$:
\begin{equation}
\ \forall \gamma \in \Gamma, \forall x \in \widetilde{S},
\: \phi(\gamma x)=\rho(\gamma)\phi(x),\nonumber
\end{equation}
and $\rho(\Gamma)$ is Fuchsian.

The \emph{number of vertices of the Fuchsian polyhedral  embedding} is the number of vertices of the cellulation 
of $S$.

The Fuchsian polyhedral embedding is \emph{convex} if its image is a convex polyhedral surface of $M_K^-$.
\end{definition}

We consider the Fuchsian polyhedral  embeddings up to homeomorphisms and up to global isometries:  $(S_1,\varphi_1,\rho_1)$ and $(S_2,\varphi_2,\rho_2)$ are equivalent if there exists 
a homeomorphism $h$ between $S_1$ and $S_2$ and an  isometry $I$ of $M_K^-$
such that, for a lift $\widetilde{h}$ of $h$ to $\widetilde{S}_1$ we have
\begin{equation}\label{def rel eq plongements}
\varphi_2\circ \widetilde{h} = I\circ \varphi_1.
\end{equation}
As two lifts of $h$ only differ by conjugation by elements of
$\Gamma$, using the equivariance property of the
embedding, it is easy to check that the definition of the
equivalence relation doesn't depend on the choice of the lift.

\begin{definition}
The \emph{genus of a Fuchsian group $F$ of $M_K^-$}  is the genus 
of the quotient of $O_K$  by the restriction of $F$ to it.

The \emph{genus of a Fuchsian polyhedron} $(P,F)$ is the genus of $F$.

The 
\emph{number of vertices of a Fuchsian polyhedron} $(P,F)$ is the
number of vertices of $P$ in a fundamental domain for the action of 
$F$.
\end{definition}
As $S$ is a compact surface of genus $g>1$, it can be endowed with hyperbolic metrics, and each of them provides 
a cocompact representation of $\Gamma$ in the group of orientation-preserving  isometries of $O_K$. The images of such representations are usually called \emph{Fuchsian groups} (as $O_K$ is isometric to $\mathbb{H}^2$), that explains the  terminology used. Using the orthogonal projection $p_K$ to extend the action of a Fuchsian group of the hyperplane $O_K$ to the entire future-cone of $c_K$, it is easy to check that (as done in \cite{artrealisationhyp}):

\begin{lemma}
There is a bijection between the cocompact representations of the fundamental group of $S$ in 
$\mathrm{Isom}^+(\mathbb{H}^2)$ and 
the Fuchsian groups of $M_K^-$ of genus $g$.

It follows that there is a bijection between the convex Fuchsian polyhedra of genus $g$ with $n$ vertices and the 
convex Fuchsian polyhedral embeddings with $n$ vertices of a compact surface of genus $g$.
\end{lemma}

\subsection{Fuchsian deformations.}\label{subsec:fuchsdef}

Let $(S,\phi,\rho)$ be a   convex  polyhedral Fuchsian embedding in $M_K^-$ and let $(\phi_t)_t$ be  a path of  convex polyhedral embeddings   of
$\widetilde{S}$ in  $M_K^-$, such that:
\begin{itemize}
\item[-] $\phi_0=\phi$,
\item[-] the induced metric is preserved at the first order at $t=0$,
\item[-] there are representations $\rho_t$ of $\Gamma=\pi_1 (S)$ into $\mathrm{Isom}_+^+(M_K^-)$
\end{itemize}
 such that
\begin{equation}\phi_t(\gamma x)=\rho_t(\gamma)\phi_t(x)\nonumber\end{equation}
and each $\rho_t(\Gamma)$ is Fuchsian. 

We denote by \begin{equation}Z(\phi(x)):= \frac{d}{dt}\phi_t(x)_{\vert t=0}\in \mathrm{T}_{\phi(x)}M_K^-\nonumber\end{equation} and
\begin{equation}\dot{\rho}(\gamma)(\phi(x))=\frac{d}{dt}\rho_t(\gamma)(\phi(x))_{\vert t=0} \in \mathrm{T}_{\rho(\gamma)\phi(x)}M_K^-.\nonumber\end{equation}
The vector field $Z$ has a property of \emph{equivariance} under
$\rho(\Gamma)$:
$$
Z(\rho(\gamma) \phi(x))=\dot{\rho}(\gamma)(\phi(x))+d\rho(\gamma).Z(\phi(x)).
$$
This can be written
\begin{equation}\label{killequiv}
Z(\rho(\gamma)\phi( x))=d\rho(\gamma).(d\rho(\gamma)^{-1}\dot{\rho}(\gamma)(\phi(x))+Z(\phi(x)))
\end{equation}
and $d\rho(\gamma)^{-1}\dot{\rho}(\gamma)$ is a Killing field of  $M_K^-$, because it is the derivative of a path in the group of isometries of  $M_K^-$ (we must multiply by  $d\rho(\gamma)^{-1}$, because
$\dot{\rho}(\gamma)$ is not a vector field). We denote this Killing field by $\vec{\rho}(\gamma)$. Equation (\ref{killequiv}) can be written, if $y=\phi(x)$,
\begin{equation}\label{equation deformation fuchsienne}
Z(\rho(\gamma) y)=d\rho(\gamma).(\vec{\rho}(\gamma) + Z)(y).
\end{equation}

 An \emph{infinitesimal isometric deformation} of a
polyhedral surface consists of 
\begin{itemize}
\item a triangulation of the polyhedral surface given by a triangulation of each face, such that no new vertex arises,
\item a Killing field on each face of the
triangulation such that  two Killing fields on two adjacent triangles are equal on the common edge. 
\end{itemize}

An infinitesimal isometric deformation is called \emph{trivial}
if it is the restriction to the polyhedral surface of a global Killing
field. If all the infinitesimal isometric deformations of a
polyhedral surface are trivial, then the polyhedral surface is said to be \emph{infinitesimally rigid}.

A \emph{Fuchsian deformation} is an infinitesimal isometric deformation  $Z$ on a Fuchsian polyhedron which satisfies Equation (\ref{equation deformation fuchsienne}), where $\vec{\rho}(\gamma)$ is a \emph{Fuchsian Killing field}, that is a Killing field of the hyperbolic plane $O_K$ extended to the future-cone of $c_K$ along the geodesics orthogonal to $O_K$. More precisely, for a point  $x\in \Omega_K$, let  $d$ be the
distance between $x$ and $p_K(x)$. We denote by $p(d)$
the orthogonal projection onto $O_K$ of the 
surface $S_d$ which is at constant distance $d$ from $O_K$ (passing
through $x$). Note that every such surface $S_d$ is orthogonal to the geodesics orthogonal to $O_K$, \emph{i.e.} containing $c_K$. Then the Killing field $L$ at $p_K(x)$ 
is extended as $dp(d)^{-1}(L)$ at the point $x$.
In other words, a Fuchsian Killing field of $M_K^-$ is a  Killing field of $M_K^-$ which restriction to each surface $S_d$  is a Killing field of $S_d$. 

A Fuchsian polyhedron  is \emph{Fuchsian infinitesimally  rigid}  if all its Fuchsian deformations are trivial (\emph{i.e.} are restriction to the Fuchsian polyhedron of Killing fields of $M_K^-$).
We want to prove
\begin{Alphatheorem}\label{theoreme de rigidite lorentz}

\begin{itemize}
\item[1)] Convex Fuchsian polyhedra in Minkowski  space are Fuchsian  infinitesimally rigid;
\item[2)] Convex Fuchsian polyhedra in de Sitter space are Fuchsian  infinitesimally rigid;
\item[3)] Convex Fuchsian polyhedra in anti-de Sitter space are Fuchsian  infinitesimally rigid.
\end{itemize}
\end{Alphatheorem}

Part 1) of Theorem~\ref{theoreme de rigidite lorentz} is done in \cite{Schpolygone} (see \cite[Thm B]{iskhakov} for a partial result).  We will deduce the two other cases by using infinitesimal Pogorelov maps defined in the next Subsection. The same method was used in the smooth case \cite{SchLab}. The key point will be that 
these infinitesimal Pogorelov maps send Fuchsian deformations to Fuchsian deformations.

Note that it is also possible to deduce Fuchsian  infinitesimally rigidity in the hyperbolic case in the same way. Another 
proof in the hyperbolic case is given in \cite{artrealisationhyp}. And it is possible to go from the hyperbolic result to the Minkowski result as in the converse way. It may also be possible to adapt the method used in  \cite{artrealisationhyp} to the de Sitter and anti-de Sitter case, but with more work as the fact to be in Lorentzian manifolds doesn't allow to do exactly the same thing as in hyperbolic case.

Note that it would be also possible to deduce de Sitter case from hyperbolic case using an argument used in \cite{CompHyp} for another type of equivariant polyhedra. The idea is   to go from hyperbolic polyhedra to de Sitter polyhedra  by projective transformations.

\subsection{Infinitesimal Pogorelov maps.}

The following construction is an adaptation of a map invented by Pogorelov 
\cite{Pogo}, which allows to transport infinitesimal deformation problems in 
a constant curvature space to infinitesimal deformation problems in a flat space,  see for example
\cite{SchLab,rousset}.  In these references, infinitesimal Pogorelov maps are introduced as the derivative of a ``global'' Pogorelov map. A more direct (equivalent) definition  was introduced in \cite{Schconvex}. We give an analogous definition here.

The \emph{radial direction} at a point $x\in\Omega_K$ is given by the derivative at $x$ of the time-like geodesic
from $c_K$ to $x$. A \emph{lateral direction} is a direction orthogonal to the radial one. Note that the radial directions are orthogonal to the surfaces of $\Omega_K$ at constant distance from $c_K$, and the lateral directions are then tangent to these surfaces.
Moreover the projective maps $\varphi_K$ sends radial directions to radial directions and lateral directions to lateral directions.

The \emph{infinitesimal Pogorelov map} from $\Omega_K$, $K=-1,1$, to Minkowski space,  is the map from $T\Omega_K$ to $T\Omega_0$ sending $(x,Z)$ to 
$(\varphi_K(x),\Phi_K(Z))$, where $\Phi_K(Z)$ is defined such that the radial component of $\Phi_K(Z)(\varphi_K(x))$  has same direction and
same norm as $Z_r(x)$, the radial component of $Z$, and the lateral component  of
$\Phi_K(Z)(\varphi_K(x))$ is $d_x\varphi_K (Z_l)$, where $Z_l$ is the lateral component of $Z$. In other words, if $R$ is the radial direction of $M^-_0$ (recall the conventions in Subsection~\ref{subsec:convention}):
$$\Phi_K(Z)(\varphi_K(x)):= d_x\varphi_K(Z_l)+ \norm{Z_r}_{\Omega_K} R(\varphi_K(x)).$$

The infinitesimal Pogorelov map has the following remarkable property:
\begin{lemma}[Fundamental property of the infinitesimal Pogorelov map]\label{pogorelov}
Let $Z$ be 
a vector field on $\Omega_K$, $K\in\{-1,1\}$, then $Z$ is a Killing field if and only if  $\Phi_K(Z)$  is a Killing field of Minkowski space.
\end{lemma}
As an infinitesimal isometric deformation of a polyhedral surface is the data of a Killing field on each triangle of a triangulation, this lemma says that
the image of an infinitesimal isometric deformation of a polyhedral surface $P$  by the infinitesimal Pogorelov map is an  infinitesimal isometric deformation of the image of $P$ by the projective map. And one is trivial when the other is.

The proof of Lemma~\ref{pogorelov} lies on the two following lemmas. We denote by $\ol{S}_{t}$ the surface of the Minkowski space at constant distance $t$ from $c_0$ and contained inside its future cone $\Omega_0$. We denote by $\ol{\I}_{t}, \ol{\II}_{t}$ respectively the induced metric on $\ol{S}_{t}$ and its second fundamental form. 

\begin{lemma}\label{spheresit2}
Let $S(\nu)$ be the surface in $M_1^-$ at constant distance $\nu$ from $c_1$ and contained in $\Omega_1$. 
Then if $t=\tanh(\nu)$ we have
$\varphi_1(S_{\nu})=\ol{S}_{t}$. Moreover if $\I_{\nu}$ and 
$\II_{\nu}$ are  respectively the induced metric on $S_{\nu}$ and its second fundamental form, then
\begin{equation*}
\I_{\nu}= \cosh^{2}(\nu)\ol{\I}_{t}, \:
\II_{\nu}=\cosh^{2}(\nu)\ol{\II}_{t}.
\end{equation*}
\end{lemma}
\begin{proof}
Let  $x\in S(\nu)$. Up to an isometry preserving $c_1$ we can consider that the time-like geodesic 
between $c_1$ and $x$ is given by the plane spanned by the coordinates  $x_1$ and $x_4$. Then 
 $x_{1}=\cosh(\nu)$ and
\begin{eqnarray*}
\ && x_{2}^{2}+x_{3}^{2}-x_{4}^{2}=1-\cosh(\nu)^{2}=-\sinh^{2}(\nu), \\
\  && \norm{\left(\frac{x_{2}}{x_{1}},
\frac{x_{3}}{x_{1}},\frac{x_{4}}{x_{1}}\right)}^2_{M^-_0}=\frac{1}{x_{1}^2}
\norm{(x_{4},x_{2},x_{3})}^2_{M^-_0}=-\frac{\sinh^2(\nu)}{\cosh^2(\nu)}
=-\tanh^2(\nu).
\end{eqnarray*}
If we denote by $\mathrm{can}_{\mathbb{H}^2}$ the induced metric on $O_1$, it follows that $\displaystyle
\ol{\I}_{t}=t^2\mathrm{can}_{\mathbb{H}^2}=\tanh^2(\nu)\mathrm{can}_{\mathbb{H}^2}$, and
it is easy to see that 
the induced metric on $S_\nu$ is 
  $$\I_{\nu}=\sinh^2(\nu)\mathrm{can}_{\mathbb{H}^2}.$$ 

The sectional curvature of $S_{\nu}$ is
$-1/\sinh^2(\nu)$, then its Gaussian curvature is  $\coth^2(\nu)$,
and as it is a totally umbilic surface we get
$$\II_{\nu}=\coth(\nu)\I_{\nu}=\cosh(\nu)\sinh(\nu)\mathrm{can}_{\mathbb{H}^2}.$$
In the same manner
$\ol{\II}_t=t\mathrm{can}_{\mathbb{H}^2}=\tanh(\nu)\mathrm{can}_{\mathbb{H}^2}$.
\end{proof}

\begin{lemma}\label{sphereads}
Let $S(\theta)$ be the surface in $M_{-1}^-$ at constant distance $\theta$ from $c_{-1}$ and contained in $\Omega_{-1}$. 
Then if $t=\tan(\theta)$ we have
$\varphi_{-1}(S_{\theta})=\ol{S}_{t}$. Moreover if $\I_{\theta}$ and 
$\II_{\theta}$ are  respectively the induced metric on $S_{\theta}$ and its second fundamental form, then
\begin{equation*}
\I_{\theta}= \cos^{2}(\theta)\ol{\I}_{t}, \:
\II_{\theta}=\cos^{2}(\theta)\ol{\II}_{t}.
\end{equation*}
\end{lemma}
\begin{proof}
Let  $x\in S(\theta)$. Up to an isometry preserving $c_{-1}$ we can consider that the time-like geodesic 
between $c_{-1}$ and $x$ is given by the plane spanned by the coordinates  $x_1$ and $x_4$. Then 
 $x_4=\cos(\theta)$ and  
\begin{eqnarray*}
\ && x_{1}^{2}+x_{2}^{2}-x_{3}^{2}=-1+\cos(\theta)^{2}=\sin^{2}(\theta), \\
\  && \norm{\left(\frac{x_{1}}{x_{4}},\frac{x_{2}}{x_{4}},
\frac{x_{3}}{x_{4}}\right)}^2_{M^-_0}=\frac{1}{x_{4}^2}
\norm{(x_{1},x_{2},x_{3})}^2_{M^-_0}=\frac{\sin^2(\theta)}{\cos^2(\theta)}
=\tan^2(\theta).
 \end{eqnarray*}
 If we denote by $\mathrm{can}_{\mathbb{H}^2}$ the induced metric on $O_{-1}$, it follows that $\ol{\I}_{t}=t^2\mathrm{can}_{\mathbb{H}^2}=\tan^2(\theta)\mathrm{can}_{\mathbb{H}^2}$, and 
 it is easy to see that 
the induced metric on $S_\theta$ is 
  $\I_{\theta}=\sin^2(\theta)\mathrm{can}_{\mathbb{H}^2}$. 

The sectional curvature of $S_{\theta}$ is
$-1/\sin^2(\theta)$, then its Gaussian curvature is $\cot^2(\theta)$, and as it is a totally umbilic surface we get $$
\II_{\theta}=\cot(\theta)\I_{\theta}=\cos(\theta)\sin(\theta)\mathrm{can}_{\mathbb{H}^2}.$$
In the same manner $
\ol{\II}_{t}=t\mathrm{can}_{\mathbb{H}^2}=\tan(\theta)\mathrm{can}_{\mathbb{H}^2}$,
\end{proof}

With this it is possible to prove that, if $g_K$ is the metric on $\Omega_K$ and $\mathcal{L}$ is the Lie derivative, $X,Y$ tangent vectors and $X',Y'$ their images by the differential of the projection, then

$$(\mathcal{L}_vg_1)(X,Y) = \cosh^{2}(\nu) (\mathcal{L}_{\Phi_1(v)}g_0)(X',Y') $$
in the de Sitter case, and in the anti-de Sitter case:
$$(\mathcal{L}_vg_{-1})(X,Y) = \cos^{2}(\theta) (\mathcal{L}_{\Phi_{-1}(v)}g_0)(X',Y'). $$

This proves Lemma~\ref{pogorelov} because a Killing field $Z$ for the metric $g$ can be defined by $(\mathcal{L}_Zg)(x,y)=0$ for all $x,y$. The way to go from lemmas \ref{spheresit2} and \ref{sphereads} to the two equations above is word by word the same as used in \cite{Schconvex,schnonconvex}. These references deal with infinitesimal Pogorelov map 
from respectively hyperbolic and de Sitter space to Euclidean space. The only difference in our case (except the constant arising in the equations above) is that the target space is the Minkowski space and not the Euclidean one, but this changes nothing, as the proof uses only definitions of Levi-Civita connection, Lie derivative and second fundamental form, which are exactly the same in Riemannian or Lorentzian geometry.

\textit{Remarks.} 
The fundamental property of  infinitesimal Pogorelov maps is related to the Darboux--Sauer Theorem, which says that the infinitesimal rigidity (for surfaces in the Euclidean space) is a property which remains true under projective transformations.
In the Riemannian case these results are contained inside a general statement  saying that each times that there is a map from a  manifold to a flat manifold sending geodesics to geodesics then there exists a (unique) map of the associated tangent bundles sending Killing fields to Killing fields \cite{volovkpogo} (this result should be checked for pseudo-Riemannian context). Concerning the polyhedral surfaces, there exists a more geometric way to define the infinitesimal Pogorelov maps \cite{white,ivanrigid}. \\

\textit{Proof of parts 2) and 3) of Theorem~\ref{theoreme de rigidite lorentz}.}
Now to prove parts 2) and 3) of Theorem~\ref{theoreme de rigidite lorentz} from its part 1) it suffices 
to show that:
\begin{lemma}\label{def fuchs sur def fuchs}
The  infinitesimal Pogorelov maps $\Phi_1$ and $\Phi_{-1}$ send Fuchsian deformations to Fuchsian deformations.
\end{lemma}
We first note that
\begin{lemma}
The  infinitesimal Pogorelov maps $\Phi_1$ and $\Phi_{-1}$  send Fuchsian Killing fields to Fuchsian Killing fields.
\end{lemma}
\begin{proof}
By Lemma~\ref{pogorelov} we know that these maps send Killing fields to Killing fields. Moreover, they send radial component to radial component, and a Fuchsian Killing field is characterised by the fact that it has no radial component.
\end{proof}
\begin{proof}[Proof of Lemma~\ref{def fuchs sur def fuchs}]
 A Fuchsian  deformation $Z$ verifies
\begin{equation}
Z(\phi( \gamma x))=d\rho(\gamma)(\vec{\rho}(\gamma)+Z)(\phi(x)),\nonumber
\end{equation}
and, on one hand the infinitesimal Pogorelov maps sends  Fuchsian Killing fields to Fuchsian  Killing fields, and on the other hand, 
 considering the radial component or  applying the projection onto the lateral component are linear operations. 
These arguments and Lemma~\ref{proj envoie fuchs sur fuchs} suffice to prove this lemma, but here are the details.

We denote by $r$ the radial direction of both $M_K^-$, $K=-1,1$, by $R$ the radial direction of $M^-_0$, by $\varphi$ the projective map $\varphi_K$ and by $\Phi$  the infinitesimal Pogorelov map $\Phi_{K}$, for $K=1$ or $K=-1$.
 Recall that the proof of Lemma~\ref{proj envoie fuchs sur fuchs} gives the existence of  a morphism $G$ between $\mathrm{Isom}_{c_K}(M^-_K)$ and  $\mathrm{Isom}_{c_0}(M^-_0)$ such that
\begin{equation}
\varphi (\rho(\gamma)(\phi(x)))=G (\rho(\gamma))(\varphi(\phi (x))).\nonumber
\end{equation}

We start from

$$ \Phi(Z)(\varphi\circ \phi (\gamma x) )=d\varphi(Z_l)(\phi(\gamma x))+\norm{Z_r(\phi(\gamma x))}R(\varphi(\phi(\gamma x))).$$

We first examine the first term of the right member of the equation above:
\begin{eqnarray*}
\ d\varphi(Z_l)(\phi(\gamma x))&=&d\varphi d\rho(\gamma)(\vec{\rho}_l(\gamma)+Z_l)(\phi(x))\\
\ &=&dG( \rho(\gamma))d\varphi(\vec{\rho}_l(\gamma)+Z_l)(\varphi \circ \phi(x)).\\
\end{eqnarray*}

Afterwards we examine the second term of the right member:
\begin{eqnarray*}
\ \norm{Z_r(\phi(\gamma x))}R(\varphi(\phi(\gamma x)))&=&\norm{d\rho(\gamma)(Z_r)(\phi(x))}R(\varphi(\rho(\gamma)\phi(x))) \\
\ &=&\norm{Z_r(\phi(x))}R(\varphi(\rho(\gamma)\phi(x))) \\
\ &=&\norm{\norm{Z_r}r(\phi(x))}R(G(\rho(\gamma))(\varphi\circ\phi(x)))\\
\ &=& \norm{Z_r}R(G(\rho(\gamma))(\varphi\circ\phi(x)))\\
\ &=& dG(\rho(\gamma))(\norm{Z_r}R(\varphi\circ\phi(x))).\\
\end{eqnarray*}

And at the end, using both computations above, we have what we wanted:
\begin{eqnarray*}
\  \Phi(Z)(\varphi\circ \phi (\gamma x) ) &=&dG( \rho(\gamma))(d\varphi(\vec{\rho}_l(\gamma))+d\varphi(Z_l)+\norm{Z_r}R)(\varphi \circ \phi(x)) \\
\ &=&dG(\rho(\gamma))(\Phi(\vec{\rho}(\gamma))+\Phi(Z))(\varphi \circ \phi(x)).
\end{eqnarray*}

Note that we used the fact that the Fuchsian isometries fix the point from which the radial direction is defined 
to write:
$$
R(G(\rho(\gamma))(\varphi\circ\phi(x)))=dG(\rho(\gamma))R(\varphi\circ\phi(x)).
$$

\end{proof}

We proved that a Fuchsian deformation in $M^-_K$ is sent to a Fuchsian deformation of $M_0^-$, $K=-1$ or $1$. Actually the infinitesimal Pogorelov maps $\Phi_K$ have inverses, which are easy to define:  they send the lateral component to its image by the inverse of the projective map and they send the radial component to a radial vector having the same norm.  These inverses 
send Fuchsian deformations of the Minkowski space  to  Fuchsian deformations of $M_K^-$, $K\in\{-1,1\}$, the proof is exactly the same as the proof of the lemma above. But the inverse of the projective map from de Sitter space to Minkowski space can send space-like polyhedral surfaces to non-space-like ones, see below. \\

\textit{Appendix: hyperbolic-de Sitter Fuchsian polyhedra.} One can define similarly an infinitesimal Pogorelov map from the hyperbolic space to the Minkowski  space. Recall that $\varphi_1$ is also a projective map from $\mathbb{H}^3_+$ to $M^-_0$ (Subsection~\ref{modelesprojectifs}).
 Actually the infinitesimal Pogorelov map is also the same as for de Sitter space: the radial directions is defined as the directions orthogonal to the surfaces at constant distance of $P_{\mathbb{H}^2}$.

We call \emph{HS-space} the union of $\mathbb{H}^3_+$ and $\Omega_1$. It is easily described in the Klein projective model, see Figure~\ref{spheredesitter}. In this picture radial directions of both spaces are the vertical ones. A polyhedral surface of the HS-space is simply a Euclidean polyhedral surface drawn in this model. The analog of Lemma~\ref{pogorelov} is proved similarly to the other cases, using the analog of Lemma~\ref{spheresit2} and Lemma~\ref{sphereads}:

\begin{lemma}
Let $S(\mu)$ be the surface in $\mathbb{H}^3$ at constant distance $\mu$ from $P_{\mathbb{H}^2}$ and contained in $\mathbb{H}^3_+$. 
Then if $t=\cotanh(\mu)$ we have
$\varphi_{1}(S_{\mu})=\ol{S}_{t}$. Moreover if $\I_{\mu}$ and 
$\II_{\mu}$ are  respectively the induced metric on $S_{\mu}$ and its second fundamental form, then
\begin{equation*}
\I_{\mu}= \sinh^{2}(\mu)\ol{\I}_{t}, \:
\II_{\mu}=\sinh^{2}(\mu)\ol{\II}_{t}.
\end{equation*}
\end{lemma}
\begin{proof}
Let  $x\in S(\mu)$. We consider  the  geodesic 
between $P_{\mathbb{H}^2}$ and $x$  given by the plane spanned by the coordinates  $x_1$ and $x_4$. Then 
 $x_1=\sinh(\mu)$ and
\begin{eqnarray*}
\ && x_{2}^{2}+x_{3}^{2}-x_{4}^{2}=-1-\sinh^{2}(\mu)=-\cosh^{2}(\mu), \\
\  && \norm{\left(\frac{x_{2}}{x_{1}},\frac{x_{3}}{x_{1}},
\frac{x_{4}}{x_{1}}\right)}^2_{M^-_0}=\frac{1}{x_{1}^2}
\norm{(x_{2},x_{3},x_{4})}^2_{M^-_0}=\frac{\sinh^2(\mu)}{\cosh^2(\mu)}
=\cotanh^2(\mu).
 \end{eqnarray*}
 If we denote by $\mathrm{can}_{\mathbb{H}^2}$ the induced metric on $P_{\mathbb{H}^2}$, it follows that $\ol{\I}_{t}=t^2\mathrm{can}_{\mathbb{H}^2}=\cotanh^2(\theta)\mathrm{can}_{\mathbb{H}^2}$, and  it is easy to see that 
the induced metric on $S(\mu)$ is 
  $\I_{\mu}=\cosh^2(\mu)\mathrm{can}_{\mathbb{H}^2}$. 

The sectional curvature of $S_{\mu}$ is
$-1/\cosh^2(\mu)$, then its Gaussian curvature is $\tanh^2(\mu)$, and as it is a totally umbilic surface we get $$
\II_{\mu}=\tanh(\mu)\I_{\mu}=\sinh(\mu)\cosh(\mu)\mathrm{can}_{\mathbb{H}^2}.$$
In the same manner $
\ol{\II}_{t}=t\mathrm{can}_{\mathbb{H}^2}=\cotanh(\mu)\mathrm{can}_{\mathbb{H}^2}$,
\end{proof}

We can even define an infinitesimal Pogorelov map from HS-space to Minkowski space by combining the ones from de Sitter and from hyperbolic space (a Killing field of one of these spaces extends uniquely to the other are they all are restrictions of Killing fields of Minkowski space of dimension $4$). The proof of Lemma~\ref{def fuchs sur def fuchs}
also extends, as the radial directions are the sames. We can then define a Fuchsian polyhedron of HS-space as a polyhedral surface of HS-space invariant under the action of a Fuchsian group fixing $c_1$ or, that is equivalent, acting on $P_{\mathbb{H}^2}$. Then part 1) of Theorem~\ref{theoreme de rigidite lorentz} gives the Fuchsian infinitesimal rigidity of the convex Fuchsian polyhedra in HS-space.

In the definition of the infinitesimal Pogorelov map from HS-space to Minkowski space introduced here, the radial direction is orthogonal to the surfaces on which the vertices of the Fuchsian polyhedra are lying. 
It is the reason why it is well adapted to the study of Fuchsian polyhedra. 
For the infinitesimal Pogorelov map from HS-space to Euclidean space, the radial direction is defined to be orthogonal to spheres, it is why it is well-adapted to the study of ``spherical'' polyhedra (``usual'' polyhedra with a finite number of vertices) \cite{schhilb,schlorentz}. 

%%%%%%%%%%%%%%%%%%%%%%%%%%%%%%
%%%%%%%%%%%%%%%%%%%%%%%%%%%%%%
%%%%%%%%%%%%%%%%%%%%%%%%%%%%%%%%

\section{Realisation of metrics}

\subsection{Sets of Fuchsian polyhedra.}\label{Set of Fuchsian polyhedra}

We denote by $\mathcal{P}_K(g,n)$ the set of convex Fuchsian
polyhedral embeddings of a compact surface $S$ of genus $g>1$ with $n$ vertices   in $M_K^-$,
modulo isotopies of  $S$ fixing the vertices of the
cellulation and modulo the isometries of $M_K^-$.
More precisely, the equivalence relation given by the isotopies is written as Equation~(\ref{def rel eq plongements}), with the difference that $\tilde{h}$ is a lift of  a homeomorphism $h$ of $S$ isotopic to the identity, such that if $h_t$ is the isotopy (\emph{i.e.} $t\in[0,1]$, 
$h_0=h$ and  $h_1=id$), then $h_t$ fixes the vertices of the cellulation for all $t$.

\begin{lemma}\label{lem:exFuchs}
 The set $\mathcal{P}_K(g,n)$ in non-empty for all $n>0$, $g>1$, $K\in\{-1,0,1\}$.
\end{lemma}
\begin{proof}

For any $g>1$, $K\in\{-1,0,1\}$, it is sufficient to prove that $\mathcal{P}_K(g,1)$ is non-empty. In this case, starting from such an element, it suffices to take a point on a face and consider its orbit under the action of the Fuchsian group. If we push all the points obtained in this way out of the polyhedron sufficiently slightly, the resulting polyhedron is still invariant by construction, and convex and space-like, as both conditions are open. We have constructed an element of $\mathcal{P}_K(g,2)$, and we can repeat the procedure.
 
Let $F$ be a Fuchsian group of genus $g$. If $K\in\{-1,0\}$, let $P$ be the boundary of the convex hull 
of the orbit for the action of $F$ of a point on $O_K$. By Lemma~\ref{fuchsien dans conelum est espace}  $(P,F)$ is an element of $\mathcal{P}_K(g,1)$.
It is not hard to construct a convex Fuchsian polyhedron of the hyperbolic space with one vertex, such that a fundamental domain for the action of the Fuchsian group is given by one face. Then  its dual provides an element of $\mathcal{P}_1(g,1)$.
\end{proof}

Let $(\varphi_1,\rho_1)$ and $(\varphi_2,\rho_2)$ be two convex Fuchsian
polyhedral embeddings describing the same element of $\mathcal{P}_K(g,n)$.  As  $h$ is homotopic to the identity,  $\forall x\in\widetilde{S}, \forall \gamma\in\Gamma$, we check from (\ref{def rel eq plongements}) 
that $ \rho_2(\gamma)(I(\varphi_1(x)))=I(\rho_1(\gamma)((\varphi_1 (x)))).$
If two orientation preserving and time orientation preserving isometries of  $M_K^-$ are equal on an open set of a totally geodesic surface
(a face of the Fuchsian polyhedron), then they are equal, it follows that for
all $\gamma\in\Gamma$, $\rho_2(\gamma)=I\circ \rho_1
(\gamma)\circ I^{-1}$. As $\rho_1$ and $\rho_2$ are also
representations of  $\Gamma$ in $\mathrm{PSL}(2,\mathbb{R})$, it follows that  they describe the same element of the Teichm\"uller space of $S$.

Moreover, it also clear from (\ref{def rel eq plongements}) that the projections onto $O_K$ of the vertices of two polyhedral surfaces given by equivalent embeddings are the same (up to a global isometry): there is a natural map   $\mathcal{S}_K$ from    $\mathcal{P}_K(g,n)$ to $\mathrm{T}_g(n)$, the Teichm\"uller space of $S$ with $n$ marked points.
From the proof of Lemma~\ref{lem:exFuchs} it is easy to see that from any Fuchsian representation and any
$n$ points on $O_K$ (in a fundamental domain), we can build a convex Fuchsian
polyhedron with $n$ vertices:  $\mathcal{S}_K$  is surjective.

To recall an element of $\mathcal{P}_K(g,n)$ from the data of its projection onto $O_K$, it remains to know 
the distance of the vertices from $c_K$. This is because a convex polyhedral surface is entirely determined by the position of its vertices (as it is the boundary of the convex hull of its vertices). Such a distance is called the \emph{height} of a vertex.

\begin{lemma}
Take $[h]\in \mathrm{T}_g(n)$. For $K\in\{-1,0\}$, $\mathcal{S}_K^{-1}([h])$ is diffeomorphic to  the open unit ball of $\mathbb{R}^{n}$. 
\end{lemma}
\begin{proof}

 We will prove that $\mathcal{S}_K^{-1}([h])$ is a contractible open subset of $(\mathbb{R}_+)^n$. We fix a fundamental domain for the action of $\rho(\Gamma)$  
on $O_K$, with $n$ marked points, which will be the projection of the vertices of the polyhedra along the geodesics from 
$c_K$ (\emph{i.e.} we fix $[h]$). We have to find the possible heights of the vertices for the resulting invariant polyhedral surface to be convex, that means that any vertex is outside of the convex hull of the other vertices. This is sufficient by Lemma~\ref{fuchsien dans conelum est espace}.

\textit{In  the anti-de Sitter space.} We consider the second Minkowski projective model of the anti-de Sitter space, where $c_{-1}$ is a point at infinity and $O_{-1}$ is the unit disc in the horizontal plane, see Subsection~\ref{modelesprojectifs}.  In this model, the convex Fuchsian polyhedra are polyhedral convex caps with vertices accumulating to the circle bounding $O_{-1}$. 
In this case we know that the set of possible Euclidean distances from $O_{-1}$ (\emph{i.e.} the horizontal plane) of the $n$ vertices of a fundamental domain is a (non-empty) contractible open subset of $(\mathbb{R}_+)^n$ \cite[Lemma 9]{artrealisationhyp}.

If $(a,b,z)$ are the Euclidean coordinates  of a  vertex, with $z$ the Euclidean distance form the horizontal plane $O_{-1}$, it is not hard to check that the anti-de Sitter distance from $c_{-1}$ of this vertex is  $$\operatorname{cot^{-1}} ( \frac{z}{\sqrt{1-a^2-b^2}}),$$ where $a$ and $b$ are fixed by hypothesis. The  set of heights is diffeomorphic to the set of Euclidean distance from $O_{-1}$: it is a contractible set.

\textit{In  the Minkowski space.}
 The proof for the anti-de Sitter space above says that the set of possible anti-de Sitter distances between $c_{-1}$ and the vertices of the polyhedral surfaces for the polyhedral surface to be convex form a contractible set. And this is always true in the Minkowski projective model, for which $c_{-1}$ corresponds to the origin. If $t$ is such an anti-de Sitter distance a direct computation shows that the corresponding Minkowski distance is $\tan(t)$: the set of possible Minkowski heights is also a contractible set.
 
\end{proof}

Recall from Subsection~\ref{Teichspace} that to each element of the Teichm\"uller space $\mathrm{T}_g$
is associated a canonical hyperbolic polygon (the Z-V-C coordinates).
A (small) open set of $\mathrm{T}_g(n)$ is parametrised by a (small) deformation
of a canonical polygon in $O_K$ and a displacement
of the marked points inside this polygon. With  fixed heights for
the vertices, a small displacement of a  convex Fuchsian
polyhedron (corresponding to a path in  $\mathrm{T}_g(n)$), is always
convex (the convexity is a property preserved by a little
displacement of the vertices) and Fuchsian (by construction).

It follows that we can endow $\mathcal{P}_K(g,n)$ with
the topology which makes it a fiber space based on $\mathrm{T}_g(n)$, with
fibers homeomorphic to a connected open subset of $\mathbb{R}^n$, and
as $\mathrm{T}_g(n)$ is a contractible manifold of dimension 
$(6g-6+2n)$:
\begin{proposition}\label{imm var hyp2}
For $K\in\{-1,0\}$, the space $\mathcal{P}_K(g,n)$ is a contractible  manifold of
dimension $(6g-6+3n)$.

For $K=1$, $\mathcal{P}_K(g,n)$ is  a connected  manifold of
dimension $(6g-6+3n)$.
\end{proposition}
\begin{proof}
It remains to prove the assertion for the de Sitter case. The description of the manifold structure is the same than for the other topologies: $\mathcal{P}_K(g,n)$ can be parametrised by little deformations of the canonical polygon, a little displacement of the marked points inside it and a little variation of their heights. This space is also path-connected: Teichm\"uller space is path-connected and vertices can be deplaced along a path in Teichm\"uller space. Some faces can meet the boundary at infinity, but one can always decrease all the heights by the same amount to avoid it.
\end{proof}

\begin{definition}
A \emph{(generalised) triangulation} of a compact  surface $S$
is a decomposition of $S$  by images by homeomorphisms of triangles of the Euclidean space, with possible
identification of the edges or the vertices, such that the interiors of the faces  (resp. of the edges) are disjoint.
\end{definition}
This definition allows triangulations of the surface with only one
or two vertices. For example, take a canonical polygon and take a vertex of this polygon. 
Join it with the other vertices of the polygon. By identifying
the edges of the polygon, we have a triangulation of the resulting
surface with only one vertex.

A simple Euler characteristic argument gives that, if $e$ is the number of edges of a triangulation, $g$ the genus of the surface and $n$ the number of vertices:
\begin{equation}
e=6g-6+3n.\nonumber
\end{equation}

A subdivision of each  faces of a convex Fuchsian polyhedron $(P,F)\in\mathcal{P}_K(g,n)$ in
triangles (such that the resulting triangulation has no more
vertices than the cellulation of the polyhedron, and is invariant
under the action of $F$) gives a triangulation denoted by $T$. Let $U$ be a small neighbourhood of $P$ in $\mathcal{P}_K(g,n)$. For each $T$ 
we introduce $\mathcal{P}^T$, which contains polyhedral surfaces  having same vertices than the elements of $U$ together with a triangulation of same combinatoric as $T$ (elements of $\mathcal{P}^T$ are Fuchsian but not necessarily convex).

We also introduce  a map $\mathrm{EP}(T)$  which sends each element of 
$\mathcal{P}^T$  to the
square  of the lengths  of  the edges of $T$ in a
fundamental domain of the Fuchsian group. As $T$ provides a triangulation of the surface $S$, the map $\mathrm{EP}(T)$ has its values in   $\mathbb{R}^{6g-6+3n}$.

\begin{lemma}\label{lem:inj loc}
 The map  
$\mathrm{EP}(T)$ is a local homeomorphism around $P$.
\end{lemma}
\begin{proof}
 The map $\mathrm{EP}(T)$ associates to the $n$ vertices $x_1,\ldots,x_n$
a set of  $(6g-g+3n)$ real numbers among all the $\mathrm{d}_{M_K^-}(fx_i,gx_j)^2$, $i,j=1,\ldots,n$, with $f,g$ elements of the Fuchsian group. It is in particular a  $C^1$  map. By the local inverse Theorem, Theorem~\ref{theoreme de rigidite lorentz}  says exactly that $\mathrm{EP}(T)$ is a local homeomorphism.
\end{proof}

\subsection{Sets of metrics.}\label{Set of metrics}

By standard methods involving Voronoi regions and Delaunay cellulations, 
it is known 
\cite{Thurart1,TroTri,Rivintriangulation}
 that
for each constant curvature metric with conical singularities on $S$ with constant sign singular curvature there exists a geodesic triangulation  such that the vertices of the triangulation are exactly the singular points. This allows us to see such a metric as a gluing of (geodesic)  triangles.
Actually, we don't need this result, because in the following we could consider only the metrics given by the induced metric 
on convex Fuchsian polyhedra. In this case, the geodesic triangulation of the metric is  given by a triangulation of the faces of the 
polyhedral surface.

We introduce the following spaces of metrics:
\begin{itemize}
\item  $\mathcal{M}(g,n)$ the set of Riemannian metrics on a compact surface $S$ of genus $g$ minus $n$ points. It is endowed with the following $C^k$ topology:
two metrics are close if their coefficients until those of their $k$th derivative in any local chart are close (we
don't care which $k > 2$);
\item $\widetilde{\mathrm{Cone}_K^{-}}(g,n)\subset \mathcal{M}(g,n)$ the space of metrics of curvature  $K$ on $S$ with $n$ conical singularities of negative singular curvature, seen as Riemannian metrics after removing the singular points; 
\item $\mathrm{Cone}_K^{-}(g,n)$ the quotient of   $\widetilde{\mathrm{Cone}_K^{-}}(g,n)$  by the isotopies of $S$ minus  $n$ marked points;
\item $\widetilde{M}^T$ --- where $T$ is a geodesic triangulation of an element of  $\widetilde{\mathrm{Cone}_K^{-}}(g,n)$ --- the space of metrics belonging to $\widetilde{\mathrm{Cone}_K^{-}}(g,n)$ which admit a geodesic triangulation homotopic to  $T$;
\item  $\mathrm{Conf}(g,n)$ the space of conformal structures on  $S$ with $n$ marked points.
\end{itemize}

We denote by $\widetilde{\mathrm{EM}}(T)$ the map from $\widetilde{M}^T$ to $\mathbb{R}^{6g-6+3n}$ which associates to each element of 
$\widetilde{M}^T$ the square of the lengths of the edges of the triangulation. 
The (square of) the distance between two marked points of $S$ is a continuous map from   $\mathcal{M}(g,n)$ to $\mathbb{R}$. 
Around a point of $\widetilde{M}^T$, $\widetilde{\mathrm{EM}}(T)$ takes its values in an open set of  $\mathbb{R}^{6g-6+3n}$: if we modify slightly the lengths of the  $(6g-6+3n)$ edges, the metric that we will obtain will always be in  $\widetilde{M}^T$, because the conditions defining a totally geodesic triangle and the ones on the values of the cone-angles  are open conditions.

\begin{lemma}\label{topologiehyp}
The space $\mathrm{Cone}_{-1}^{-}(g,n)$ is a contractible manifold of dimension $(6g-6+3n)$.
\end{lemma}
\begin{proof}

Theorems of Picard--Mc Owen--Troyanov \cite{owen}\cite[Theorem A]{Troyanovarticle2} say that, if $g>1$,
there is a bijection between  $\widetilde{\mathrm{Cone}_{-1}^{-}}(g,n)$ and $\mathrm{Conf}(g,n)\times A_n$,   where $A_n$ is a contractible sets of $\mathbb{R}^{n}$  given by Gauss--Bonnet conditions ($A_n$ parametrises the values of the cone-angles). 

As the Teichm\"uller space $\mathrm{T}_g(n)$ is the quotient of  $\mathrm{Conf}(g,n)$  by the isotopies of $S$ minus its marked points,
 $\mathrm{Cone}_{-1}^{-}(g,n)$ is in bijection with $\mathrm{T}_g(n)\times A_n$, and this last space is contractible. With the help of this bijection, we endow $\mathrm{Cone}_{-1}^{-}(g,n)$ with the topology of $\mathrm{T}_g(n)\times A_n$.
\end{proof}

\begin{lemma}\label{topologieplate}
The space $\mathrm{Cone}_{0}^{-}(g,n)$ is a contractible manifold of dimension $(6g-6+3n)$.
\end{lemma}
\begin{proof}

A Theorem of  Troyanov \cite{Troyanovarticle1,Troyanovarticle2} says that, if $g>1$,
there is a bijection between    $\widetilde{\mathrm{Cone}_{0}^{-}}(g,n)$ up to the homotheties and $\mathrm{Conf}(g,n)\times B_n$, where  $B_n$ is a contractible sets of $\mathbb{R}^{n-1}$  given by Gauss--Bonnet conditions ($B_n$ parametrises the values of the cone-angles).
We endow $\mathrm{Cone}_{0}^{-}(g,n)$ with the topology which makes it a fiber space based on $\mathrm{T}_g(n)\times B_n$ with fiber $\mathbb{R}_+$.
 \end{proof}

\begin{lemma}
For $K\in\{-1,0\}$, let $M^T$ be the quotient of  
$\widetilde{M}^T$ by the isotopies, and  $\mathrm{EM}(T)$ the map induced  by $\widetilde{\mathrm{EM}}(T)$. 
Then $\mathrm{EM}(T)$  is a local homeomorphism from $M^T\subset\mathrm{Cone}_K^{-}(g,n)$ to its image in $\mathbb{R}^{6g-6+3n}$.
\end{lemma}
\begin{proof} 
For the topology given by the one of the space of metrics,   $\widetilde{\mathrm{EM}}(T)$ is a continuous map on  $\widetilde{M}^T\subset \widetilde{\mathrm{Cone}_K^{-}}(g,n)$.
For the cases $K\in\{-1,0\}$, we check that this property is always true for the topologies given by Lemma~\ref{topologiehyp} and Lemma~\ref{topologieplate}.
   
  Let  $i_T$ be the canonical inclusion of $\widetilde{M}^T$ (endowed with the topology induced by the one of  $\mathcal{M}(g,n)$) in $\widetilde{\mathrm{Cone}_{-1}^{-}}(g,n)$. 
  For the hyperbolic case, the composition of  $i_T$ with the projection onto   $\mathrm{Conf}(g,n)$ 
   is the map which associates to each metric its conformal structure, 
   this is a continuous map as by definition $\mathrm{Conf}(g,n)$ is the quotient of 
$\mathcal{M}(g,n)$ by the set of real-values functions on $S$ minus its marked points.
Moreover, the composition of $i_T$ with the projection onto   $A_n$ is obviously continuous.
It follows that $i_T$ is continuous and injective: it is a local homeomorphism.  

It is exactly the same for the flat case, because if we fix a point $m$ in $\widetilde{M}^T$, around it $\widetilde{\mathrm{Cone}_{0}^{-}}(g,n)$ can be written as $\mathrm{Conf}(g,n)\times B_n\times \mathbb{R}_+$. 
To conclude it remains to note that, up to the isotopies of the surface, the map  $\widetilde{\mathrm{EM}}(T)$ becomes an injective map  $\mathrm{EM}(T)$ from $M^T$ 
to $\mathbb{R}^{6g-6+3n}$. This map takes its values in an open set of   $\mathbb{R}^{6g-6+3n}$ and  the dimension of  $\mathrm{Cone}_K^{-}(g,n)$ is $(6g-6+3n)$, 
 \end{proof}

Things are not so simple for the spherical metrics, because if there exists a result of existence of the metrics (under a Gauss--Bonnet condition), the uniqueness is not known \cite{Troyanovarticle2} --- actually, there exists uniqueness results for some particular cases, see  \cite{fengluo,umehara,eremenko}.
  For this reason, we can endow  $\mathrm{Cone}_1^-(g,n)$ only with the topology given by the one of $\mathcal{M}(g,n)$.
  We denote by $\mathrm{Cone}_1^{-,>2\pi}(g,n)$ the subset of $\mathrm{Cone}_1^-(g,n)$ containing the metrics with closed contractible geodesics of lengths $>2\pi$.
  
\begin{lemma}\label{ens met cone spher loc manifold}
The space $\mathrm{Cone}_1^{-,>2\pi}(g,n)$ is locally a manifold of dimension   $(6g-6+3n)$.
\end{lemma}
\begin{proof}
With the help of a triangulation $T$, using $\widetilde{\mathrm{EM}}(T)$ as  above, we get that  the space $\mathrm{Cone}_1^{-}(g,n)$ is locally a manifold of dimension $(6g-6+3n)$. 
Moreover, the condition on the lengths of the closed contractible geodesics is an open condition \cite[Theorem 6.3, Lemma 9.9]{RivHod}.
\end{proof}

%%%%%%%%%%%%%%%%%%%%%%%%%%%%%%%%%%%%%%%%%%%%%%%%%%%%%%%%%%%
\subsection{Final steps.}\label{final steps}

 We denote by $\mathcal{I}_K(g,n)$ the map  which associates to each element $(P,F)$ of 
$\mathcal{P}_K(g,n)$ the induced metric on $P/F$, which is known to be an element of $\mathrm{Cone}^-_K(g,n)$. In the spherical case, the induced metric lies more precisely in $\mathrm{Cone}_1^{-,>2\pi}(g,n)$. An element of $\mathrm{Cone}^-_K(g,n)$ lying in the image of $\mathcal{I}_K(g,n)$ is said to be a \emph{realisable} metric.
\begin{lemma}
 The map
$\mathcal{I}_K(g,n)$ is continuous and locally injective.
\end{lemma}
 \begin{proof}
 The map
$\mathcal{I}_K(g,n)$ is obviously continuous.
Suppose that $\mathcal{I}_K(g,n)$ is not locally injective. That means that there exists  $(P',F')$ as close as we want from $(P,F)$ such that the induced metric on both are isometric. If they are sufficiently close, we can endow them with the same triangulation $T$. Then they belong to the same space $\mathcal{P}^T$ and their image under $\mathrm{EM}(T)$ have to be different by Lemma~\ref{lem:inj loc}. On the other hand it is obvious that
$$ \mathrm{EM}(T)\circ \mathcal{I}_{K}(g,n)\circ \mathrm{EP}(T)^{-1}= id,
$$
that gives a contradiction with the fact that the induced metrics are isometric. 
\end{proof}

The map $\mathcal{I}_{K}(g,n)$ is continuous and locally injective, hence it is a local homeomorphism by the invariance of domain Theorem. 
In the next section we will show that $\mathcal{I}_K(g,n)$ is proper, and as $\mathrm{Cone}^-_{K}(g,n)$ is connected it follows that $\mathcal{I}_{K}(g,n)$ is surjective, therefore a covering map. As  $\mathrm{Cone}^-_{K}(g,n)$ is simply connected and  $\mathcal{P}_K(g,n)$ is connected, 
it follows  that 
$\mathcal{I}_{K}(g,n)$ is a homeomorphism.

Let $\mathrm{Mod}(n)$ be the quotient of the group of the homeomorphisms of $S$
minus $n$ points by its subgroup of isotopies.
The homeomorphism $\mathcal{I}_{K}(g,n)$ gives a bijection between
 the quotient of $\mathcal{P}_{K}(g,n)$ by $\mathrm{Mod}(n)$ and the quotient of $\mathrm{Cone}^-_{K}(g,n)$ by $\mathrm{Mod}(n)$ for $K\in\{-1,0\}$, and this is exactly the statement of
  parts 2) and 3) of Theorem~\ref{realisation lorentz}.

 Now it remains to prove that  $\mathcal{I}_1(g,n)$ is 
a homeomorphism between $\mathcal{P}_{1}(g,n)$ and $\mathrm{Cone}_1^{-,>2\pi}(g,n)$. But we don't know anything about the connectedness  of $\mathrm{Cone}_1^{-,>2\pi}(g,n)$: the conclusion is less straightforward than for Minkowski and anti-de Sitter spaces.
In \cite{rivinthese,RivHod} there is a result on a kind of ``connectedness'' for $\mathrm{Cone}_1^{-,>2\pi}(0;n)$, using the  connectedness of a space of smooth metrics, and  it is noted in \cite{Schpoly} that the genus doesn't intervene in the proof. The only difference is that in our case we must consider the metrics up to isotopies, that changes nothing.

\begin{proposition}\label{sorte de connexite metriques dans de sitter}
Each metric $m_1\in \mathrm{Cone}_1^{-,>2\pi}(g,n)$ can be joined to  a metric 
$m_0:=\mathcal{I}_1(g,n)(P)$, for a $P\in \mathcal{P}_{1}(g,n)$, by a continuous path $(m_t)_t$, with $m_t\in  \mathrm{Cone}_1^{-,>2\pi}(g,N)$,  $N\geq n$, $t\in ]0,1[$, and such that $m_t$ is realisable for $t$ near $0$.
\end{proposition}
\begin{proof}[Sketch of the proof]
For a suitable neighbourhood of the cone points of $m_0$ and $m_1$, it is possible to (continuously) smooth each cone point \cite[9.2]{RivHod}  to obtain continuous paths $(\overline{m}_t)_t$, $t\in[0,t_1]$ and $(\overline{m}_{t'})_{t'}$, $t'\in[t_2,1]$, where $\overline{m}_{t_1}$ and $\overline{m}_{t_2}$ are smooth metrics with curvature $K\leq 1$ and lengths of contractible geodesics $L>2\pi$  (obviously, $\overline{m}_0=m_0$ and $\overline{m}_1=m_1$). The space of such metrics is path-connected (that is proved using standard arguments \cite{RivHod,schconvlor,SchLab,Schconvex}).
It comes that $m_0$ and $m_1$ can be joined by a continuous path of (smooth or with conical singularities) metrics such that $K \leq 1$ and $L>2\pi$. 

Now take a geodesic cellulation of $m_0$ such that the cone points are the vertices, and subdivide each cell with as many (geodesic) triangles as necessary to each triangles to have a diameter strictly less than a given constant $\delta$. We denote by $N$ the number of vertices resulting of such a triangulation $T_0$. The deformation $(\overline{m}_t)_t$ gives a continuous family $T_t$ of geodesic triangulations, and $T_1$ is a geodesic triangulation of $m_1$. Afterward we replace each triangle by a spherical triangle with the same edge length, and this gives us the announced path $m_t$ between $m_0$ and $m_1$ (this new path can be taken very close to $(\overline{m}_t)_t$, such that its cone angles remain $>2\pi$ and the lengths of its closed contractible geodesics remain $> 2\pi$).
 
It remains to prove that for $t$ sufficiently small, $m_t$ is realisable. The triangulation of $m_0$ gives a triangulation of $P$, and each $m_t$, $t\in[0,\epsilon]$ is obtained by pushing outward each vertex (of the triangulation) contained inside a face of $P$. The way to push each vertex is given by the change of the length of the edges of the triangulation. This technique is also used in \cite{Aleks}.
\end{proof}

In the next section, we will show that $\mathcal{I}_1(g,n)$ is proper, and thus:

\begin{corollary}
The map $\mathcal{I}_1(g,n)$ is surjective.
\end{corollary}
\begin{proof}
With the same notations than above,  we already know that $m_t$ is realisable for $t\in[0,\epsilon[$. By properness of $\mathcal{I}_1(g,N)$,  $m_t$ is realisable for $t\in[0,\epsilon]$. By local injectivity and the fact that  $\mathrm{Cone}_1^{-,>2\pi}(g,n)$ is locally an open manifold, the invariance of domain Theorem gives that the map $\mathcal{I}_1(g,N)$ is open: $m_t$ is realisable for $t\in[0,\epsilon'[$, with $\epsilon'>\epsilon$, and so on. 
At the end,  $m_{t}$ is realisable for $t\in [0,1[$, and again by properness of $\mathcal{I}_1(g,N)$, $m_1$ is realisable.
\end{proof}

The following way to conclude that $\mathcal{I}_1(g,n)$ is a homeomorphism is proposed in  \cite{Schpoly}. We know that $\mathcal{I}_1(g,n)$ is a covering on the entire $\mathrm{Cone}_1^{-,>2\pi}(g,n)$ and that $\mathcal{P}_{1}(g,n)$ is connected. To conclude that $\mathcal{I}_1(g,n)$ is a homeomorphism, it remains to check that each fiber contains only one element. This is equivalent to prove that the covering of a loop is a loop, using a kind of ``simple connectedness'' of $\mathrm{Cone}_1^{-,>2\pi}(g,n)$, and this is given by  a straightforward adaptation of  Proposition~\ref{sorte de connexite metriques dans de sitter}:
\begin{proposition}
For each  $c : \mathbb{S}^1 \rightarrow \mathrm{Cone}_1^{-,>2\pi}(g,n)$ there exists a disc $D\subset  \mathrm{Cone}_1^{-,>2\pi}(g,N)$,  $N\geq n$, $t\in ]0,1[$,  such that $\partial D = c(\mathbb{S}^1)$.
\end{proposition}
Note that we know now that all the metrics involved in this lemma are realisable.
\begin{proof}[Sketch of the proof]
The proof is step by step the same as for Proposition \ref{sorte de connexite metriques dans de sitter}, using the fact that the space of smooth metrics with curvature $\leq 1$ and lengths of contractible geodesics $>2\pi$ is simply connected --- that is proved using standard arguments \cite{RivHod,schconvlor,SchLab,Schconvex}.
\end{proof}

%%%%%%%%%%%%%%%%%%%%%%%%%%%%%%%%%
\section{Properness}\label{properness}

We will use the following characterisation of a proper map: 
$\mathcal{I}_K(g,n)$ is proper if, for each sequence $(P_k)_k$ in $\mathcal{P}_K(g,n)$ such that 
the sequence $(g_k)_k$ converges in $\mathrm{Cone}^-_K(g,n)$ (with $g_k:=\mathcal{I}_K(g,n)(P_k)$) to $g_{\infty}\in \mathrm{Cone}^-_K(g,n) $, then 
$(P_k)_k$ converges in $\mathcal{P}_K(g,n)$ (may be up to the extraction of a sub-sequence). For $K=1$, we consider $\mathrm{Cone}_1^{-,>2\pi}(g,n)$ instead of $\mathrm{Cone}^-_1(g,n)$.

We denote by  $\mathrm{d}_k$ the restriction  to $P_k$ of the distance from $c_K$. We always denote by $p_K$ the orthogonal projection onto $O_K$, and $(\varphi_k,\rho_k)$ is the embedding of the surface $S$ corresponding to $P_k$.
Let denote by  $\gamma_k$  a geodesic on $P_k$ given by an element of the fundamental group of $S$ or   a geodesic between 
two vertices, and by $l_k$ the length of $\gamma_k$. By convergence of the sequence of induced metrics, $l_k$ is bounded from above and below for all $k$. We denote these bounds by $l_{\mathrm{min}}\leq l_k \leq l_{\mathrm{max}}$. Note that this argument will avoid the collapsing of  singular points. 
We will suppose that the geodesics $\gamma_k(t)$ are parametrised by the arc-length, \emph{i.e.}  $g_k(\gamma_k'(t),\gamma_k'(t))=1$.

For each $M^-_K$, we call $u_k$ the restriction of a  coordinate function to $P_k$, that is
\begin{equation}\label{definitions fonctions coord}
u_k:=\left\lbrace
\begin{array}{ll}
  \frac{1}{2}(\mathrm{d}_k)^2,\: K=0; \\
 \cos (\mathrm{d}_k),\: K=-1;\\
\cosh(\mathrm{d}_k),\: K=1.
\end{array}
\right.
\end{equation}

\begin{lemma}\label{lemme pas de max sur les aretes}
For $K\in\{0,1\}$, for each $k$, for each geodesic $\gamma(t)$ on $P_k$,  $(u_k\circ \gamma)'$ has a positive jump at its singular points (which correspond to points where $\gamma(t)$ crosses an edge of $P_k$).
\end{lemma}
\begin{proof}

Consider an edge $e$ of $P_k$. We denote by $f_1$and $f_2$ its adjacent faces, and we look at a geodesic $\gamma (t)$ (for the induced metric) on $f_1\cup f_2$. We denote by $\gamma_i$ the part  of $\gamma$ which lies on $f_i$, and $t_0$ is such that $\gamma (t_0)\in e$. We denote by $\overline{\gamma}_1$ the prolongation of 
$\gamma_1$ on the plane containing the face $f_1$. Let $d$ be the distance from $c_K$. The graph of $d\circ \overline{\gamma}_1$ is  smooth, and, until $t_0$, the graph of $d\circ \gamma$ is also smooth. 

At $t_0$, $d\circ \gamma_2$ and $d\circ \overline{\gamma}_1$ have the same value, and   as $P_k$ is convex  and $c_K$ lies in the concave side of $P_k$,  from $t_0$ $d\circ \gamma_2$ is greater  than $d\circ \overline{\gamma}_1$, that means that the jump of $(d\circ \gamma)'$ at $t_0$ is positive.
As the geodesic lies on $P_k$, we can write that the jump of $(\mathrm{d}_k \circ \gamma)'$ at $t_0$ is positive, and as the functions involved in (\ref{definitions fonctions coord}) are increasing for $K\in\{0,1\}$, this is true for $(u_k\circ \gamma)'$.
\end{proof}

\begin{lemma}\label{distance bornee de en bas uniforme}
For all $k$, the distance $\mathrm{d}_k$ is uniformly bounded from below by a strictly positive constant.
\end{lemma}
\begin{proof}\cite{Schpoly}
We see a sequence of (closure of) fundamental domains on $P_k$ for the action of $\rho_k(\Gamma)$ as a sequence $(D_k)_k$ of convex isometric space-like embeddings of the disc, with $n$ singular points. Each $D_k$ must stay out of the light-cone of its vertices by Lemma~\ref{lem:outlightcone}, and inside the light-cone of $c_K$. It follows that if a vertex $x_k$ goes to $c_K$, then the $D_k$ will be in an arbitrarily neighborhood of a light-cone for $k$ sufficiently large. But this is impossible: a light-cone (without its vertex) is a smooth surface, and it cannot be approximate by polyhedral surfaces with a fixed number of vertices.
\end{proof}
\begin{lemma}\label{dilatation convergence}
If the projection of the $P_k$ onto a space-like  surface $N$ at constant distance from $c_K$ is a dilating function, then the associated sequence of representations converges.
\end{lemma}
\begin{proof}

The curvature of the induced metric on $N$ is constant and strictly negative. 
For each $k$,  $\rho_k(\Gamma)$ acts on $N$, and the quotient is isometric to a hyperbolic metric (up to a homothety) on the compact surface $S$. 
We denote by $h_k$ this hyperbolic metric on  $S$.

By hypothesis, the induced metrics $g_k$ on $P_k$  converge to $g_{\infty}$. For $k$ sufficiently large, there exists a constant $c'$ such that $g_k\geq \frac{1}{c'}g_{\infty}$. As the surface is compact, there exists a constant $c$ such that
$\frac{1}{c'}g_{\infty}\geq \frac{1}{c}h_{0}$. And as the projection is dilating, we have, if $\mathrm{L}_g(\gamma)$ is the length of the geodesic corresponding to $\gamma\in\pi_1(S)$ for the metric $g$:
\begin{equation}
\mathrm{L}_{h_k}(\gamma) \geq \mathrm{L}_{g_k}(\gamma) \geq \frac{1}{c}\mathrm{L}_{h_0}(\gamma),\nonumber
\end{equation}
and  Lemma \ref{lemme teichmuller minor geod}  leads to the conclusion.
\end{proof}

\begin{lemma}\label{rep et haut suffisent}
Up to extract a subsequence, if the sequence of representation $(\rho_k)_k$ converges and the height of at least one vertex is bounded from above, then the sequence $(P_k)_k=(\varphi_k,\rho_k)_k$ converges in $\mathcal{P}_K(g,n)$.
\end{lemma}
\begin{proof}
We first check that the positions of the vertices (\emph{i.e.} $(\varphi_k)_k$) converge. As the Fuchsian embeddings are defined up to global isometries and as we suppose that they lie inside the future-cone of $c_K$, up to compose on the left by a sequence of  isometries of the future-cone of $c_K$, we can consider that there exists a vertex $x_k\in P_k$ which always remains on the same geodesic from $c_K$. As the representations converge, for $k$ sufficiently large, all the vertices (in a fundamental domain) are contained inside a cone, built with the images of $x_k$  under the action of generators of the fundamental group of the surface.
 
 We can consider that the heights   which are bounded from above are those of $x_k$. For each $k$, consider the convex hull $C_k$ of $x_k$ together with  the orbits of $x_k$ under the action of the Fuchsian group. By hypothesis, these convex hulls converge to a convex polyhedral surface $C$. As  $P_k$ is convex and  as $c_K$ lies in the concave side of $P_k$, all the vertices of $P_k$ must lie  in the same side of $C$ than $c_K$. It follows that the heights of all the vertices are bounded from above, and also  from below by Lemma~\ref{distance bornee de en bas uniforme}. It follows that the vertices (for a fundamental domain)  are contained inside a truncated cone, that is a compact domain.

It follows that the sequence $(P_k)_k$ converges  to a  polyhedron $P$ invariant under the action of a Fuchsian group of genus $g$, which must be convex with $n$ vertices as the $P_k$ are convex and as sequence of the induced metrics converges. It remains to check that $P$ is space-like, as some faces could become light-like. In Minkowski and anti-de Sitter cases, $P$ is space-like by  Lemma~\ref{fuchsien dans conelum est espace}. In the de Sitter case, it is known that if one face becomes light-like, then a geodesic of length $2\pi$ bounding a hemisphere appears in the limit metric \cite{RivHod}, that is impossible as the metrics are supposed to converge in $\mathrm{Cone}_1^{-,>2\pi}(g,n)$.
\end{proof}

\subsection{Properness in the anti-de Sitter space.}

In  $\Omega_{-1}$, the anti-de Sitter  metric can be written 
$\sin^2(t)\mathrm{can}_{\mathbb{H}^2}-dt^2$, where $t$ is the distance to $c_{-1}$ and $\mathrm{can}_{\mathbb{H}^2}$ the  hyperbolic metric.
  In the projective model for which $c_{-1}$ is sent to infinity, all the $P_k$ lie above $O_{-1}$ (Lemma~\ref{region des pol fuch ads}), this  means that the projection onto $O_{-1}$ is dilating and by Lemma~\ref{dilatation convergence}, the sequence of representations associated to the $P_k$ converges.
Moreover, in this model, the heights of the vertices are bounded, as all the $P_k$ lie below the  surface realising the minimum of the distance to $c_{-1}$ and above the horizontal plane. Lemma~\ref{rep et haut suffisent} leads to the conclusion.

\subsection{Properness in the de Sitter space.}\label{Properness in the de Sitter space}

Almost of this part  was done in \cite{Schpoly}.

\begin{lemma}\label{preuvecompactdesit}
The sequence of the representations associated to $(P_k)_k$ converges (up to extract  a subsequence).
\end{lemma}
\begin{proof}

Let $\gamma$ be an element of the fundamental group of $S$ as in Lemma~\ref{lemme teichmuller major geod}.
At this $\gamma$ corresponds a minimising geodesic $\gamma_k(t)$ on $P_k$ between a point $\varphi_k(x)\in P_k$ and $\varphi_k(\gamma x)\in P_k$. We denote by $L_k$ the length of the projection of $\gamma_k(t)$ onto $O_1$. If we prove that $L_k$ is bounded from above for all $k$, then Lemma~\ref{lemme teichmuller major geod} will lead to the conclusion.

We denote by $g_{k}$ the induced metric on $P_k$, which can be written:
\begin{equation}
g_k=\sinh^2 (\mathrm{d}_k)\mathrm{can}_{\mathbb{H}^2}- d\mathrm{d}_k^2\nonumber
\end{equation}
that leads to
\begin{equation}
 g_k=(u_k^2-1)\mathrm{can}_{\mathbb{H}^2}-\frac{du_k^2}{u_k^2-1},\nonumber
\end{equation}
it follows that we can compute:
\begin{eqnarray*}
\ L_k&=&\int_0^{l_k} \sqrt{\mathrm{can}_{\mathbb{H}^2}(dp_1(\gamma_k'(t)),dp_1(\gamma_k'(t)))}dt\nonumber \\ 
\ &=& \int_0^{l_k} \sqrt{\frac{g_k(\gamma_k'(t),\gamma_k'(t))}{u_k^2(t)-1}+\frac{du_k^2(\gamma'(t))}{(u_k^2(t)-1)^2}} dt\\
\ &=& \int_0^{l_k} \sqrt{\frac{1}{u_k^2(t)-1}+\frac{u_k'^2(t)}{(u_k^2(t)-1)^2}} dt \nonumber \\
\ &=&  \int_0^{l_k} \sqrt{\frac{1}{u_k^2(t)-1}+(\operatorname{cotanh^{-1}}(u_k(t))')^2} dt \\
\ &\leq & \int_0^{l_{\mathrm{max}}} \sqrt{\frac{1}{u_k^2(t)-1}+(\operatorname{cotanh^{-1}}(u_k(t))')^2} dt \\
\ &\leq & \int_0^{l_{\mathrm{max}}}  \sqrt{\frac{1}{u^2_k(t)-1}} dt+\int_0^{l_{\mathrm{max}}} \vert (\operatorname{cotanh^{-1}}(u_k(t)))'\vert dt \\
\ &\leq &   \frac{l_{\mathrm{max}}}{\sqrt{u^2_0-1}}+\int_0^{l_{\mathrm{max}}} \vert (\operatorname{cotanh^{-1}}(u_k(t)))'\vert dt 
\end{eqnarray*}
(we have used the fact that  $u_k$ is bounded from  below by $u_0>\cosh(0)=1$,  Lemma~\ref{distance bornee de en bas uniforme}).

It remains to check that the variation of $ \operatorname{cotanh^{-1}}(u_k)$ over $[0,l_{\mathrm{max}}]$ is bounded from above by a constant which does not depend on $k$. For this, we can decompose $[0,l_{\mathrm{max}}]$ into a finite number of subsets of the form $[x,y]$, where $x$ is a local minimum (of $u_k$) and $y$ a local maximum, which immediately follows $x$ in the list of local extrema, and into a finite number of subsets of the form $[y,x]$, where $y$ is a local maximum and $x$ a local minimum, which immediately follows $y$ in the list of local extrema. 

First we consider a subset of the kind $[x,y]$, where $x$ is a local minimum and $y$ a local maximum, which immediately follows $x$ in the list of local extrema. We want to study the variation
\begin{equation}\label{premiere variation}
\int_x^y \vert (\operatorname{cotanh^{-1}}(u_k(s)))'\vert ds=\vert \operatorname{cotanh^{-1}}(u_k(x))-\operatorname{cotanh^{-1}}(u_k(y))\vert.
\end{equation}

There exists a brutal overestimation which is:
\begin{equation}
\vert \operatorname{cotanh^{-1}}(u_k(x))-\operatorname{cotanh^{-1}}(u_k(y)) \vert \leq \operatorname{cotanh^{-1}}(u_0),\nonumber
\end{equation}
but it is not satisfying: as the number of subsets in the decomposition of $[0,l_{\mathrm{max}}]$ actually depends on $k$, the bound may become very large together with $k$. We will use the above bound only in the case $\vert y-x\vert \geq \pi /4$, and we will compute another bound in the other case (the term $\pi / 4$ has no particular role, the important thing is that the other case verifies $\vert y-x\vert < \pi /2$).

As $y$ is a local maximum, together with Lemma~\ref{lemme pas de max sur les aretes}, we have $u'_k(y)=0$. 

We also know  that, if $f$ is the restriction to a pseudo-sphere of a linear form (and $u_k$ are such functions) then $\mathrm{Hess} (f) = -f g$, where $g$ is the induced metric on the pseudo-sphere \cite[Ex. 2.65,b]{GaHuLa}. It 
 gives that $u_k''=-u_k$ on the regular points, and, again by Lemma~\ref{lemme pas de max sur les aretes}, the derivative has a positive jump at the singular points. With these facts, it is easy to check that, for $s\in[x,y]$,
\begin{equation} 
u_k(y)\geq u_k(s)\geq u_k(y) \cos(y-s),\nonumber
\end{equation}
and with this we compute
\begin{eqnarray*}
\ \vert \operatorname{cotanh^{-1}}(u_k(x))-\operatorname{cotanh^{-1}}(u_k(y)) \vert &\leq& \int_{u_k(x)}^{u_k(y)} \frac{dt}{t^2-1}\ 
\\ &\leq& \frac{u_k(y)-u_k(x)}{u_k(x)^2-1} \\
\ &\leq& \frac{u_k(x)}{u_k(x)^2-1}\left( \frac{1}{\cos(y-x)}-1 \right) \\
\ &\leq& 4 (y-x)^2\frac{u_k(x)}{u_k(x)^2-1} \
\\ &\leq& \frac{4(y-x)^2u_0}{u_0^2-1}.
\end{eqnarray*}

The bound is the same in the case where a local minimum immediately follows a local maximum, by the symmetry in  Equation~(\ref{premiere variation}).  At the end we have the wanted bound:
\begin{equation}
\int_0^{l_{\mathrm{max}}}\vert (\operatorname{cotanh^{-1}}(u_k(t)))'dt\vert
\leq \frac{4l_{\mathrm{max}}(\operatorname{cotanh^{-1}}(u_0)-1)}{\pi}+\frac{4u_0l^2_{\mathrm{max}}}{u_0^2-1}.\nonumber
 \end{equation}
\end{proof}

\begin{lemma}\label{lem:desitfaces espace}
The heights of the vertices converge.
\end{lemma}
\begin{proof}
We want to prove  that the height of no vertex of the $P_k$  goes to  infinity. It is known that in this case the sequence of induced metrics will converge to a metric having  a geodesic of length $2\pi$ \cite{RivHod,Schpoly}, that is impossible as the induced metrics are supposed to converge in $\mathrm{Cone}_1^{-,>2\pi}(g,n)$. 
\end{proof}

The conclusion follows from Lemma~\ref{rep et haut suffisent}.

\subsection{Properness in the Minkowski space.}

\begin{lemma}
The sequence of the representations associated to $(P_k)_k$ converges (up to extract  a subsequence).
\end{lemma}
\begin{proof}
We will  prove it as it had be done for Lemma~\ref{preuvecompactdesit}. We take back the same notations as given in the first lines of the proof of this lemma.

It is easy to check that the induced metric $g_k$ on $P_k$ can be written:
\begin{equation}
 g_k=u_k\mathrm{can}_{\mathbb{H}^2}-\frac{du_k^2}{u_k},\nonumber
\end{equation}
it follows that we can compute:
\begin{eqnarray}
\ L_k&=&\int_0^{l_k} \sqrt{\mathrm{can}_{\mathbb{H}^2}(dp_0(\gamma_k'(t)),dp_0(\gamma_k'(t)))}dt \nonumber \\
\ &=& \int_0^{l_k} \sqrt{\frac{g_k(\gamma_k'(t),\gamma_k'(t))}{u_k(t)}+\frac{du_k^2(\gamma_k'(t))}{u_k^2(t)}} dt \nonumber \\
\ &=& \int_0^{l_k} \sqrt{\frac{1}{u_k(t)}+\left(\frac{u_k'(t)}{u_k(t)}\right)^2} dt\nonumber \\
 \ &\leq&  \int_0^{l_{\mathrm{max}}}  \sqrt{\frac{1}{u_k(t)}+\left(\frac{u_k'(t)}{u_k(t)}\right)^2} dt\nonumber \\ 
\ &\leq & \int_0^{l_{\mathrm{max}}}  \frac{1}{\sqrt{u_k(t)}}+\left\vert\frac{u_k'(t)}{u_k(t)}\right\vert dt. \label{unemajration pour proprete de sitter fuchsien2}
\end{eqnarray}

 As  $u_k$ is bounded from  below by $u_0>0$ (Lemma~\ref{distance bornee de en bas uniforme}), we get 
\begin{equation}
L_k \leq  \frac{l_{\mathrm{max}}}{\sqrt{u_0}}+ \int_0^{l_{\mathrm{max}}} \vert (\ln(u_k(t))'\vert dt, \nonumber
\end{equation}
and it remains to check that the variation of $\ln(u_k(t))$ over $[0,l_{\mathrm{max}}]$ is bounded from above by a constant which does not depend on $k$. 

We introduce the same decomposition of $[0,l_{\mathrm{max}}]$ than for the de Sitter case:  we decompose $[0,l_{\mathrm{max}}]$ into a finite number of subsets of the form $[x,y]$, where $x$ is a local minimum  (of $u_k$) and $y$ a local maximum, which immediately follows $x$ in the list of local extrema, and into a finite number of subsets of the form $[y,x]$, where $y$ is a local maximum and $x$ a local minimum, which immediately follows $y$ in the list of local extrema. 

Without loss of generality, we suppose that $u_0>1$.
First we consider a subset of the kind $[y,x]$, where $x$ is a local minimum and $y$ a local maximum, which immediately follows $x$ in the list of local extrema. We  want to study the variation
\begin{equation}\label{premiere variation2}
\int_y^x \vert (\ln(u_k(t)))'\vert dt=\ln(u_k(y))-\ln(u_k(x)).
\end{equation}

 By Lemma
\ref{lemme pas de max sur les aretes}, $u_k'(y)=0$. Furthermore, $u_k''=-1$ on the regular points ($u_k$ is defined as (half) minus the squared norm, its Hessian is minus the bilinear form). Moreover, $u_k'$ has a positive jump at certain points.  From these facts, it is easy to check that, for $s\in [0,x-y]$: 
\begin{equation}\label{minoration hauteur mink}
u_k(y+s)\geq u_k(y)-\frac{s^2}{2},\nonumber \end{equation} 
in particular,
\begin{equation}\label{minoration hauteur mink2}
u_k(y)-u_k(x)\leq \frac{(x-y)^2}{2}\nonumber. \end{equation} 

From this we compute
\begin{equation}\label{unemajration pour proprete de sitter fuchsien}
  \ln(u_k(y))-\ln(u_k(x))  \leq  \int_{u_k(x)}^{u_k(y)} \frac{dt}{t} \leq \frac{u_k(y)-u_k(x)}{u_k(x)}\leq \frac{(x-y)^2}{2u_0}.
\end{equation}

The bound is the same in the case where a local maximum immediately follows a local minimum, by the symmetry in  Equation (\ref{premiere variation2}). It implies that
\begin{equation}
\int_0^{l_{\mathrm{max}}}\vert (\ln(u_k(t)))'\vert dt
\leq \frac{(l_{\mathrm{max}})^2}{2u_0}.\nonumber
\end{equation}
\end{proof}

\begin{lemma}
The height of at least one vertex converges.
\end{lemma}
\begin{proof}
From (\ref{unemajration pour proprete de sitter fuchsien2}) and (\ref{unemajration pour proprete de sitter fuchsien}) we can write:
$$
L_k\leq  \frac{l_{\mathrm{max}}}{\sqrt{\mathrm{min}_k}}+\frac{(l_{\mathrm{max}})^2}{2 \mathrm{min}_k},
$$ 
where $\mathrm{min}_k$ is, for each $k$, the minimum of the $u_k(x)$, where $x$ is a local minimum for the restriction of $u_k$ to $\gamma_k$ (it may be not unique). If the height of no vertex is bounded,  $\mathrm{min}_k$ will become big when $k$ is large,  and $L_k$ will be close to $0$. But   all the $L_k$ (built for each  $\gamma\in\pi_1(S)$) can't be arbitrarily close to $0$, because if it is, the area of a fundamental domain on $\mathbb{H}^2$ for the action of the Fuchsian representations will be close to $0$,  that is impossible by the Gauss--Bonnet Formula.
\end{proof}

The conclusion follows from Lemma \ref{rep et haut suffisent}.

\bibliographystyle{alpha}

\end{document}